\documentclass{siamart171218}

\usepackage{graphicx}
\usepackage{amsmath}
\usepackage{graphics}
\usepackage{amsfonts}
\usepackage{amssymb}%
\usepackage{amssymb}
\usepackage{mathrsfs}
\usepackage{amsfonts}
\usepackage{amsmath}
\usepackage{bm}
\usepackage{graphicx}
\usepackage{multirow}
\usepackage{caption2}
\usepackage{float}
\usepackage{booktabs}
\usepackage{algorithm}
\usepackage{algorithmic}
\usepackage{subfigure}

\usepackage{array}

\usepackage{amsmath}
\usepackage{amsfonts}
\usepackage{graphicx}
\usepackage{epstopdf}
\usepackage{comment}
\usepackage{subfigure}
\usepackage{algorithmic}
\ifpdf
  \DeclareGraphicsExtensions{.eps,.pdf,.png,.jpg}
\else
  \DeclareGraphicsExtensions{.eps}
\fi

\newsiamremark{remark}{Remark}
\newsiamremark{hypothesis}{Hypothesis}
\crefname{hypothesis}{Hypothesis}{Hypotheses}
\newsiamthm{claim}{Claim}

\headers{Distributed Uncertainty Quantification of Kernel Interpolation on Spheres}{S.-B. Lin,   X.Sun and D. Wang}

\title{Distributed Uncertainty Quantification of Kernel Interpolation on Spheres \thanks{The corresponding author is Di Wang. 
\funding{S. B. Lin  was partially
	supported by the  National Key R\&D Program of China (No.2020YFA0713900) and the
    Natural Science Foundation of China [Grant No  62276209].}}}

\author{Shao-Bo Lin\thanks{Center for Intelligent Decision-Making and Machine Learning, School of Management, Xi'an Jiaotong University, Xi'an 710049, China
  (\email{sblin1983@gmail.com}).}
\and Xingping Sun\thanks{ Department
	of Mathematics, Missouri State University, Springfield, MO 65897,
	USA
  (\email{XSun@MissouriState.edu}).}
  \and Di Wang\thanks{Center for Intelligent Decision-Making and Machine Learning, School of Management, Xi'an Jiaotong University, Xi'an 710049, China
  (\email{wang.di@xjtu.edu.cn}).}
  }

\usepackage{amsopn}

\newcommand{\sfgrad}[1][]{ 
	\nabla_{*}
}

\newcommand{\sfcurl}[1][]{ 
	\mathbf{L}
}

\newcommand{\Jw}[1][\alpha,\beta]{ 
w_{#1}
}

\newcommand{\imat}[1][d]{ 
    I
}

\newcommand{\Lpw}[2][\Jw]{ 
\mathbb{L}_{#2}(#1)
}

\newcommand{\InnerLGb}[2][{\Jw[r-\frac{1}{2},r-\frac{1}{2}]}]{ 
\left(#2\right)_{\Lpw[{#1}]{2}}
}

\newcommand{\Diff}[2][t]{ 
\ifthenelse{\equal{#2}{1}}{\frac{\mathrm{d}}{\mathrm{d}#1}}{
\left(\frac{\mathrm{d}}{\mathrm{d}#1}\right)^{#2}}
}

\ifpdf
\hypersetup{
  pdftitle={Distributed Kernel Interpolation on Spheres},
  pdfauthor={S.-B. Lin, X. Sun and D. Wang}
}
\fi

\begin{document}

\maketitle
\begin{abstract}
For radial basis function (RBF) kernel interpolation of scattered data, Schaback \cite{schaback1995error} in 1995 proved that the attainable approximation error and the condition number of the underlying interpolation matrix cannot be made small simultaneously. He referred to this finding as an ``uncertainty relation",  an undesirable consequence of which is that RBF kernel interpolation is susceptible to noisy data. In this paper, we propose and study a distributed interpolation method to manage and quantify the uncertainty brought on by interpolating noisy spherical data of non-negligible magnitude. We also present numerical simulation results showing that our method is practical and robust in terms of handling noisy data from challenging computing environments. 
\end{abstract}




\pagestyle{myheadings}
\thispagestyle{plain}

\section{Introduction}

Numerical methods pertaining to scattered spherical data analysis are increasingly needed and applied in a large swath of scientific fields including  geophysics \cite{king2012lower}, planetary science \cite{wieczorek1998potential}, computer graphic science \cite{tsai2006all} and  signal recovery \cite{mcewen2011novel}.
In many of the real world application problems,  data collected are of the form 
$D:=\{( x_i, y_i)\}_{i=1}^{|D|}$, in which $x_i \in \mathbb S^d$, the unit sphere embedded in the $(d+1)$-dimensional Euclidean space $\mathbb R^{d+1}$, $ y_i \in \mathbb R$, and $|D|$ denotes the cardinality of the data set $D$.
Data analysts aim to develop new and\slash or apply existing  numerical algorithms to generate an estimator $f_D$ that approximates or fits the unknown but definitive
relation between inputs and outputs. We are interested in the following  noisy data fitting model on spheres:
\begin{equation}\label{Model1:fixed}
        y_{i}=f^*(x_{i})+\varepsilon_{i},  \qquad\forall\
        i=1,\dots,|D|,
\end{equation}
where $\{\varepsilon_{i}\}_{i=1}^{|D|}$ is a set of independent zero-mean random noises satisfying $|\varepsilon_i|\leq M$ for some  $M>0$ and $f^*$ is a function from $\mathbb S^d$ to $\mathbb R$ which we use to model the relation between the input $x$ and output $y$.  

Let $\phi$ be a positive definition function \cite{hubbert2015spherical} on $\mathbb S^d$ and $(\mathcal N_\phi,\|\cdot\|_\phi)$ the native space (also called  reproducing kernel Hilbert space) associated with $\phi$.  Minimal norm kernel interpolation defined by
\begin{equation}\label{minimal-norm-interpolation}
   f_{D}:=  \arg\min\limits_{f\in\mathcal N_\phi}\|f\|_\phi,\qquad
    \mbox{s.t.} \quad 
     f(x_i)=y_{i},\quad (x_{i},y_{i})\in D, 
\end{equation}
is a classical and long-standing approach for spherical scattered data fitting, in which $f_{D}$ is explicitly written as 
\begin{equation}\label{KI}
    f_{D}=\sum_{i=1}^{|D|}a_i\phi_{x_i}, \qquad \mbox{with}\
   (a_1,\dots,a_{|D|})^T= :{\bf a}_{D}=\Phi_D^{-1}{\bf y}_D.
\end{equation} 
Here
$\phi_{x_i}$ denotes the function on $\mathbb S^d$: $ x \mapsto \phi(x_i\cdot x)$, ${\bf y}_D=(y_1,\ldots, y_{|D|})^T$, and $\Phi_D$ is the kernel matrix $(\phi(x_i\cdot x_j))_{i,j=1}^{|D|}$ in which  $x \cdot y $ is the usual dot product of $x$ and $y$. Customarily, it requires  orders of complexity  $\mathcal O(|D|^3)$ and $\mathcal O(|D|^2)$ in training and storage respectively to get such an estimator as in \eqref{minimal-norm-interpolation}. The ``meshless" feature of   kernel interpolation has rendered itself useful in practice and attracted researchers' attention on its stability in numerical implementation or the lack thereof \cite{narcowich1998stability,levesley1999norm}. Furthermore, approximation error estimates for kernel interpolation have been established separately under the circumstances of fitting clean data \cite{narcowich2002scattered,narcowich2007direct,hangelbroek2010kernel,hangelbroek2011kernel} and  noisy data \cite{hesse2017radial,hesse2021local,Feng2021radial}.

\begin{figure*}[t]
    \centering
    \subfigure[RMSE and sample size]{\includegraphics[width=6cm,height=4cm]{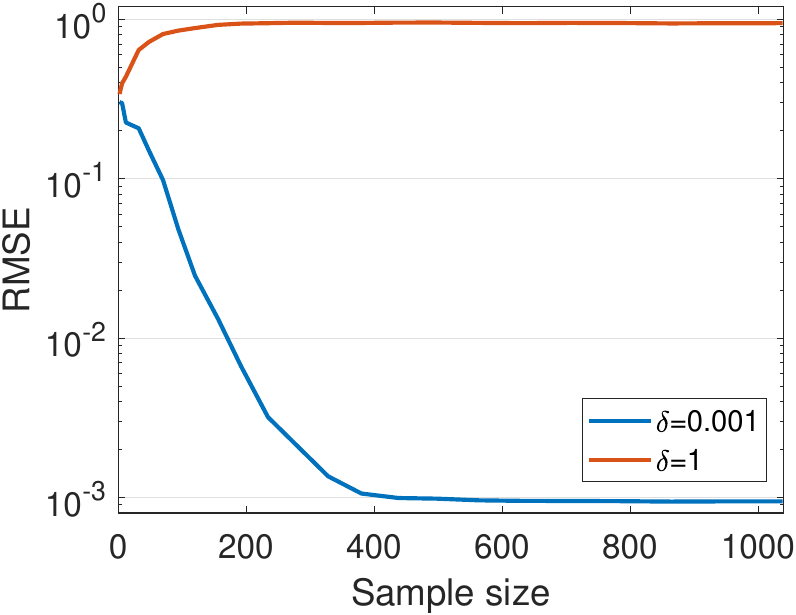}}
     \subfigure[CNKM and sample size]{\includegraphics[width=6cm,height=4cm]{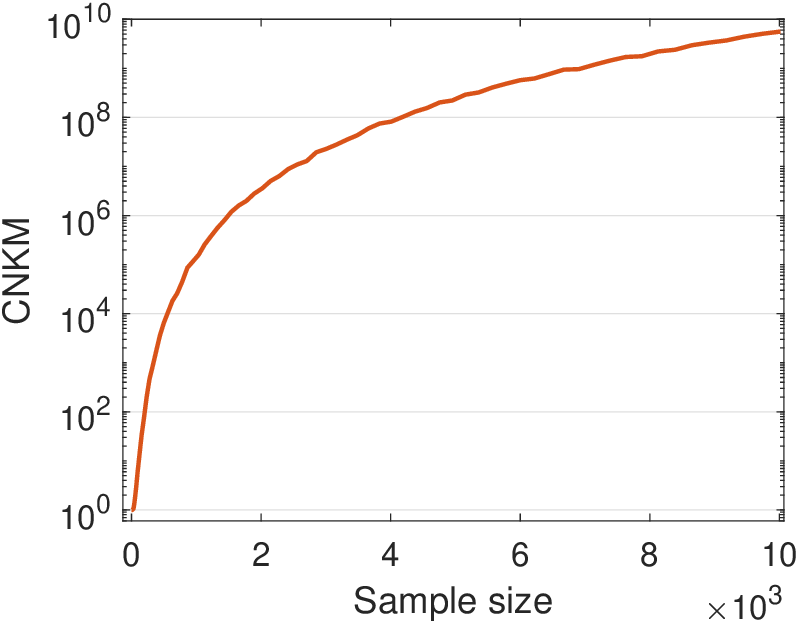}}
    \caption{Relations between training sample size and RMSE (rooted 
   mean square error)/CNKM (condition number of the kernel matrix). 
    The training sample inputs $\{{x}_i\}_{i=1}^N$ are generated by Womersley's symmetric spherical $t$-designs \cite{womersley2018efficient} on $\mathbb S^2$, and the corresponding outputs are of the  form $y_i = f(x_i)+\varepsilon_i$, where $f$ is defined by (\ref{f1}) below and $\varepsilon_i$ are independent Gaussian noise $\mathcal{N}(0, \delta^2)$. 
   We ran simulations 30 times and recorded the average RMSE  using the kernel given as in (\ref{psi_wave}) below.
   }
    \label{Fig: RMSE and condition}
\end{figure*}

Despite its strength in utility, the stability in numerical implementation of kernel interpolation has been a lingering concern. In 1995, Schaback used the term ``uncertainty relation" to describe the dichotomy between the approximation error and condition number of kernel matrix $\Phi_D$ \cite{schaback1995error,wendland2004scattered}. (We will further elaborate on this topic in the next section.) In a nutshell, 
Schaback's uncertainty relation reveals  
  that the approximation error and condition number of kernel matrix $\Phi_D$ cannot be made small simultaneously. An undesirable consequence of the uncertainty relation is that kernel interpolation is susceptible to  noisy data. Furthermore, our numerical simulations (Figure \ref{Fig: RMSE and condition}) show that increasing sample size disproportionately escalates severity of the uncertainty relation and negatively affects the performance of kernel interpolation.  

A widely adopted approach to circumvent uncertainty in kernel interpolation is the so called Tikhonov regularization which incorporates  an additional regularization parameter $\lambda>0$ to balance the approximation error and condition number of the kernel matrix $\Phi_D+\lambda |D| I$. 
Hesse et al \cite{hesse2017radial} and  Feng et al \cite{Feng2021radial} analyzed the approximation error of Tikhonov regularization fitting spherical data respectively with small and deterministic noise and large and random noise.  
It is worth noting that deterministic noise and random noise
are merely two contrasting perspectives which data analysts take  in their effort to quantify uncertainty in kernel interpolation, as are
the mathematical tools utilized to tackle them.

This paper proposes and studies a new method: distributed kernel interpolation (DKI). The basic idea of DKI is to break up judiciously the data set $D$ into $m$ blocks $D_1, \ldots, D_m,$ where $m \in \mathbb N.$ 
Based on each $D_j (1 \le j \le m)$, we build an ``interpolet" $f_{D_j}$. The final product is a global estimator $f_D$ in the format of \eqref{DKI}, which is effectively a convex linear combination of interpolets $f_{D_j}$. We carry out both 
theoretical analysis and numerical experiments in showing that DKI is a robust way of fitting massive and noisy spherical data.


The rest of the paper is organized as follows. In Section 2, we give  a detailed account of the uncertainty relation under the setting of kernel interpolation on spheres. In Section \ref{Sec:dis}, we propose the DKI method and present results of our theoretical analysis. Section \ref{sec.sampling-ineq} provides essential tools and preparations for  the proofs in Section \ref{Sec.Proofs}, highlighting an employment of a certain kind of quadrature rules of optimal orders incorporating an integral operator approach.  Section \ref{Sec.Numerical} devotes to the presentation of  four numerical simulation results. 

\section{Uncertainty Relation of Kernel Interpolation on Spheres}\label{Sec.Uncertainty}
Let $L^2(\mathbb S^d)$  be the Hilbert space with inner product
\[
\langle f, g\rangle = \int\limits_{\mathbb S^d} f(x) g(x) d \omega(x), \quad f,g \in L^2(\mathbb S^d),
\]
where $d\omega$ denotes the rotational invariant probability  measure on $\mathbb S^d$.
Let $\mathbb{H}^{d}_k$ be the subspace of $L^2(\mathbb S^d)$ consisting of all the homogeneous spherical harmonics of degree $k$. Suppose that $\{Y_{k,\ell}\}_{\ell=1}^{Z(d,k)}$ is an orthonormal basis of $\mathbb H_k^d$. Here  $Z(d,k):= \mbox{dim}\ \mathbb{H}^{d}_k$. Let $P_k^{d+1}$ denote the Gegenbauer polynomial of order $\frac{d-1}{2}$ and degree $k$ normalized so that $P_k^{d+1}(1)=1$. It is known \cite{muller2006spherical} that 
\begin{equation}\label{dimension}
 Z(d,k)=\left\{\begin{array}{ll}
  \frac{2k+d-1}{k+d-1}{{k+d-1}\choose{k}}, & \mbox{if}\ k\geq 1, \\
1, & \mbox{if}\ k=0,
\end{array}
\right.
\end{equation}
and that  
\begin{equation}\label{addition-formula}
             \sum_{\ell=1}^{Z(d,k)}Y_{k,\ell}(x)Y_{k,\ell}(x')=\frac{Z(d,k)}{\Omega_{d}}P_k^{d+1}(x\cdot x'), \quad x, x' \in \mathbb S^d,
\end{equation}
in which $\Omega_d=\frac{2\pi^{\frac{d+1}{2}}}{\Gamma(\frac{d+1}{2})}$ is the volume of $\mathbb S^d$. Equation \eqref{addition-formula} is recognized in the literature as the addition formula for homogeneous spherical harmonics. For $f \in L^2(\mathbb S^d), $ its Fourier coefficients $\hat{f}_{k,\ell}$ are defined as $
                 \hat{f}_{k,\ell}:=\int_{\mathbb
                 S^d}f(x)Y_{k,\ell}(x)d\omega(x).
$
We take the liberty of assuming that the Fourier coefficient definition above is extended in a routine manner so that it is applicable to all the tempered distributions on $\mathbb S^d.$

For a nonnegative integer $s$, denote by  $\mathcal P_s^{d}$ the set of  all   algebraic polynomials with degree at most $s$ defined on $\mathbb S^d$. Then we have that $\mathcal P_s^{d}=\bigoplus_{k=0}^s\mathbb{H}^{d}_k$ and that dim$\mathcal P_s^d = 
\sum_{k=0}^sZ(d,k)=Z(d+1,s) \sim s^d$. Here and hereafter by writing $A(p) \sim B(p)$ between the two nonnegative quantities $A(p)$ and $B(p)$ depending on a parameter $p$ other than the dimension $d$, we mean that there exist two constants $0<C_1<C_2$ depending only on $d$ such that $C_1 B(p) \le A(p)\le C_2 B(p). $
A univariate tampered distribution $\phi$ supported on $[-1,1]$ is called a spherical basis distribution (SBD) \cite{narcowich2007direct}, if
its Gegenbauer series expansion
$
           \phi(u)=\sum_{k=0}^\infty
            \hat{\phi}_k\frac{Z(d,k)}{\Omega_d} P_k^{d+1}(u)
$ has all positive coefficients, i.e., 
$0<\hat{\phi}_k< \infty$ for all $k=0,1,\ldots.$ Here the series convergence is under the topology of tampered distributions. 
If in addition  $\sum_{k=0}^\infty\hat{\phi}_k {Z(d,k)} <\infty$, then $\phi$ is called  a spherical basis function (SBF). Given an SBD $\phi$, we define
$$
               \mathcal N_\phi:=\left\{f(x)=\sum_{k=0}^\infty\sum_{\ell=1}^{Z(d,k)}\hat{f}_{k,\ell}Y_{k,\ell}(x): \|f\|_\phi^2:=
               \sum_{k=0}^\infty
         \hat{\phi}_k^{-1}\sum_{\ell=1}^{Z(d,k)}\hat{f}_{k,\ell}^2<\infty\right\} 
$$ 
as the native space of $\phi$. In particular, if $\hat{\phi}_k=1$ for all $k,$ then $\mathcal N_\phi$ coincides with $L^2(\mathbb S^d)$.

The approximation power of the kernel interpolation scheme defined as in \eqref{minimal-norm-interpolation} depends on the distribution of the  point set of inputs $\Lambda:=\{x_i\}_{i=1}^{|D|}$ on $\mathbb S^d$ and the spectrum of $\Phi_D$; see e.g. \cite{narcowich2002scattered,narcowich2007direct}. We follow what are prevalent in the scattered data approximation literature in defining three parameters that quantify the distribution of  points $x_1, \ldots, x_{|D|}$ on $\mathbb S^d$: the mesh norm $h_{\Lambda}$, the separation radius $q_{\Lambda}$, and the mesh ratio $\rho_{\Lambda}$:
\begin{align}
 h_{\Lambda}&:=\max_{ x\in\mathbb S^d}\min_{ x_{i}\in \Lambda}\mbox{dist}( x,x_i), \\
 q_{\Lambda}&:=\frac12\min_{i\neq i'}\mbox{dist}( x_{i}, x_{i'}), \\
 \rho_{\Lambda}&:=\frac{h_{\Lambda}}{q_{\Lambda}}.  
\end{align}
Here $\mbox{dist}( x,x' )$
 is the geodesic (great circle) distance between   $x$ and $x'$ on $\mathbb S^d$. If there is a constant $\tau=\tau(d)$ depending only on $d$, such that $\rho_{\Lambda}\leq \tau$, then we say that $\Lambda$  is $\tau$-quasi-uniformly distributed on $\mathbb S^d,$ or simply that $\Lambda$  is $\tau$-quasi-uniform. An elementary sphere-packing argument shows that such a $\tau(d)$ must grow exponentially with the dimension $d.$  We point out here that 
 the mesh norm is actually the Hausdorff distance between the two sets $\Lambda$ and $\mathbb S^d.$ In the sequel, we will also refer to the points in $\Lambda$ as sampling sites.

If we have in  \eqref{Model1:fixed} that $\varepsilon_i=0$ almost surely, then we refer to this as the ``noise-free model", which 
 has been extensively studied in the literature; see e.g. \cite{le2006continuous,narcowich2002scattered,narcowich2007direct,hangelbroek2010kernel,hangelbroek2011kernel}. In particular, the following result was proved in \cite[Theorem 5.5]{narcowich2007direct}.

\begin{lemma}\label{Lemma:Narcowich-interpolation}
Suppose that $\phi$ is an SBF with $\hat{\phi}_k\sim k^{-2\gamma}$ with $2\gamma>d$. Let $\alpha$ and $\beta$ be nonnegative real numbers satisfying $0\leq \beta\leq \alpha\leq 1$  and $\alpha\gamma>d/2$. Suppose that $\psi$ and $\varphi$ are SBFs such that
 \begin{eqnarray}\label{kernel-relation}
    \hat{\varphi}_k \sim \hat{\phi}_k^{\alpha},\quad \hat{\psi}_k \sim 
     \hat{\phi}_k^{\beta} 
\end{eqnarray} 
If $f^*\in \mathcal N_\varphi$, then the noise-free model in \eqref{Model1:fixed} of kernel-$\phi$ interpolation yields the error estimate
\begin{equation}\label{Error-est-noiseless}
     \|f^*-f_D\|_{\psi}\leq C \rho_\Lambda^{1-\alpha}
     h_\Lambda^{(\alpha-\beta)\gamma},
\end{equation}
 where $C$ is a constant depending only on $d, \alpha, \beta,$ and $\gamma$.   
\end{lemma}

However, a different picture emerges when one deals with noisy data. As is well-known from basic numerical linear algebra, the stability or the lack thereof for the kernel interpolation algorithm depend on the condition number of the kernel matrix $\Phi_D$. In the RBF setting, Shaback \cite{schaback1995error} established the following interesting inequality:
\begin{equation}\label{uncertainty}
   P^2(x) \left(\sigma_{\min}(\Phi^{(x)}_D)\right)^{-1} \ge 1, \quad x \in \mathbb S^d \setminus \{x_1, \ldots, x_{|D|}\}.
\end{equation}
Here $\sigma_{\min}(A)$ denotes the smallest  eigenvalue of a positive definite matrix $A,$ $\Phi^{(x)}_D$ the
$(|D|+1) \times (|D|+1)$ kernel matrix associated with the $(|D|+1) $ points  $x, x_1, \ldots, x_{|D|}$, and $P(x)$ the so called ``power function" defined by
\[
P(x):= \| \phi_x - \mathcal{P}_\Lambda (\phi_x)\|_\phi,
\]
in which $\mathcal{P}_\Lambda$ denotes the projector from $\mathcal{N}_\phi$ to its subspace $S_\Lambda:=\mbox{span}
\{\phi_{x_j}\}^{|D|}_{j=1}.$ Due to the reminiscence of \eqref{uncertainty} to the uncertainty principle in quantum mechanics, Schaback referred to it as the uncertainty relation of kernel interpolation. Three timely remarks are in order. 1. The power function $P(x)$ has a direct bearing on the rate of attainable approximation of the kernel interpolant; 2.  Schaback's uncertainty relation  carries over to the SBF setting in a routine fashion; 3. The two quantities $ \sigma_{\min}(\Phi^{(x)}_D)$ and $\sigma_{\min}(\Phi_D)$ are indistinguishably close for appropriately chosen $x \in \mathbb S^d \setminus \{x_1, \ldots, x_{|D|}\}.$

Simply put, if the rate of attainable approximation of the kernel interpolant is $R(h_\Lambda)$, where $R: (0,\infty) \rightarrow (0,\infty)$ is an increasing function and satisfies $\lim_{t \downarrow 0}R(t)=0,$ then Schaback's uncertainty relation asserts that $\left(\sigma_{\min}(\Phi_D)\right)^{-1} $ grows at least the magnitude  $R^{-1/2}(h_\Lambda)$, resulting in a noise propagation of a sizable proportion and accentuating the notion that for the type of noisy data specified in \eqref{Model1:fixed}, the approximation error of kernel interpolation of our concern based on quasi-uniformly distributed sampling sites does not decrease statistically when the sample size increases. The following result and our numerical experiments further lend support to Schaback's uncertainty relation; see Figure \ref{Fig: RMSE and condition}.

\begin{proposition}\label{Corollary:inconsistence}
Let $\Lambda$ be $\tau$-quasi-uniform for some $\tau>1$. Suppose that $\hat \phi_k\sim k^{-2\gamma}$  with $\gamma> d/2$. Let  $\{\varepsilon_i^*\}_{i=1}^{|D|}$ be a set of random variables whose supports are contained in $[-M, -\theta M] \cup [\theta M, M] $, where $0< \theta <1$ is an absolute constant.
Then for all $f \in \mathcal N_\varphi$ with $f|_\Lambda =0,$ there holds almost surely the following inequality:
\begin{equation}
     \|f_{D}-f\|_\phi\geq \tilde{C},
\end{equation}
where $f_D$ is defined by \eqref{minimal-norm-interpolation} with $y_i=f(x_i)+\varepsilon^*_i$ and
$\tilde{C}$ is a positive constant depending only on $\theta, \tau,\gamma,M,$ and $d$.
\end{proposition}

\section{Distributed Kernel Interpolation}\label{Sec:dis}
Corollary \ref{Corollary:inconsistence} shows that kernel interpolation performs poorly while confronting noisy data of non-negligible magnitude. To overcome this major drawback,  we propose and study in this section a distributed kernel interpolation  (DKI) method, which is motivated by the ``distributed learning"  in the literature \cite{zhang2015divide,lin2017distributed}. Figuratively speaking, this is a divide-and-conquer strategy for uncertainty quantification. To elaborate, we 
 describe the method in three steps. 
 \begin{itemize}
     \item 
 {\it Step 1. Data set decomposition:} For a fixed $1 \le m \le |D|$, decompose the data set $D$ into $m$ disjoint subsets $D_j:=\{(x_{i,j},y_{i,j}\}_{i=1}^{|D_j|}$ so that $D=\bigcup_{j=1}^mD_j$, $D_j\cap D_{j'}=\varnothing$ for $j\neq j'$, and the input set $\Lambda_j$ of $D_j$ is $\tau$-quasi-uniform for some $\tau>2$.

\item  {\it Step 2. Kernel interpolation on data sub-sets:} On each subset $D_j$, carry out the minimal norm interpolation as follows. \begin{equation}\label{minimal-norm-interpolation-local}
   f_{D_j}:=  \arg\min\limits_{f\in\mathcal N_\phi}\|f\|_\phi,\qquad
    \mbox{s.t.} \quad 
     f(x_{i,j})=y_{{i,j}},\quad (x_{i,j},y_{i,j})\in D_j. 
\end{equation}
We will henceforth refer to the interpolant $f_{D_j}$ as an ``interpolet". 

\item {\it Step 3. Synthesization of interpolets:} Build the DKI estimator by synthesizing interpolets via weighted averaging
\begin{equation}\label{DKI}
      \overline{f}_{D}=\sum_{j=1}^{m}\frac{|D_j|}{|D|}f_{D_j}.
\end{equation}
 \end{itemize}
 
\begin{figure}[t]
  \centering
  \begin{tabular}{m{3.1cm}<{\centering} m{0.7cm}<{\centering} m{3.1cm}<{\centering} m{0.7cm}<{\centering} m{3.1cm}<{\centering}}
  \includegraphics[width=3.1cm,height=3cm]{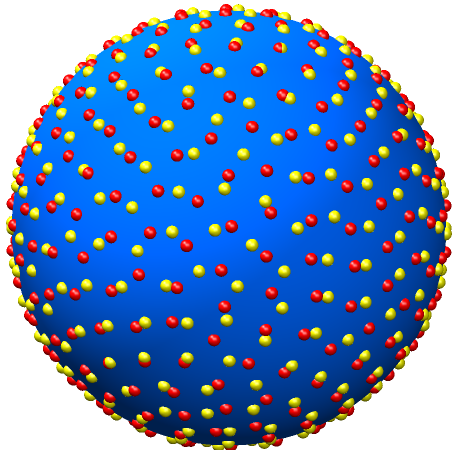}
  &\Large{$\Rightarrow$}
  &\includegraphics[width=3.1cm,height=3cm]{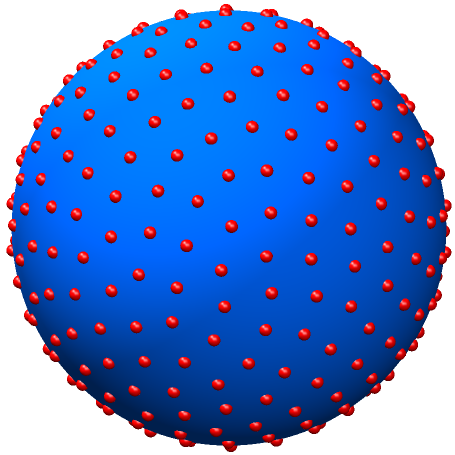}
  &\Large{$+$}
  &\includegraphics[width=3.1cm,height=3cm]{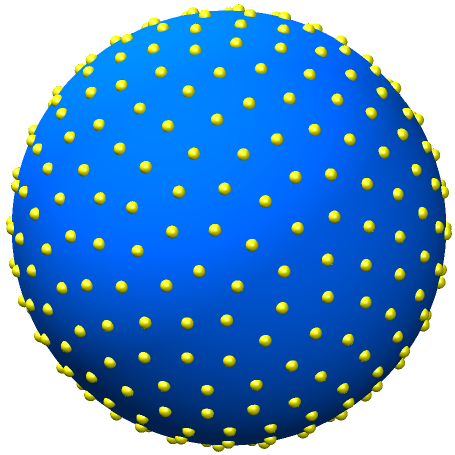}
  \end{tabular}
  \caption{The division of samples.}\label{Fig:Division}
\end{figure}

\begin{remark}
    Step 1 encompasses decomposing   the whole data set $\Lambda$ as the union of $m$ disjoint $\tau$-quasi-uniform sets.  Figure \ref{Fig:Division} exhibits such  an example. For a general $\tau$-quasi-uniform data set, we adopt a  selection-and-judgement (SAJ) strategy for our decomposition purpose, for which the algorithm entails an order $\mathcal O(|D|^2)$ of complexity in computation. We present a thorough analysis of the algorithm   in Appendix A of the Supplementary Material of the paper. 
    
\end{remark}

\begin{figure*}[t]
    \centering
    {\includegraphics[width=6cm,height=4cm]{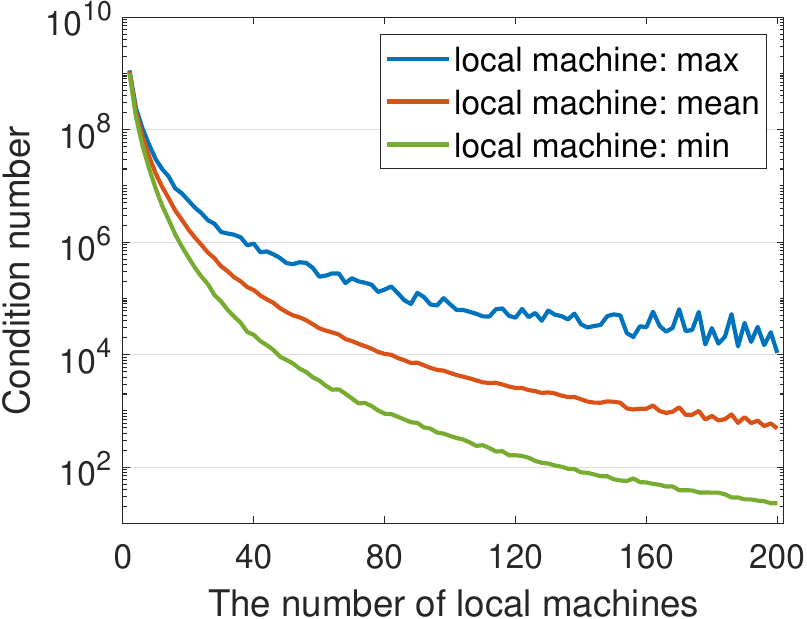}}
    \caption{Relation between the number of local machines and the average condition number of kernel matrices from local machines. The simulation setting is the same as that in Figure \ref{Fig: RMSE and condition}.}\label{Fig:dis-condition}
\end{figure*}

 DKI in \eqref{DKI} has two major advantages over the classical KI.  From the computational point of view,  DKI requires training complexity  of order $\mathcal O(\sum_{j=1}^m|D_j|^3)$, which is considerably less than  $\mathcal O(|D|^3)$ as required by KI. Furthermore,   implementing  DKI requires $m$ matrix inversions respectively for
$\Phi_{D_j},j=1,\dots,m$. For each $j$, $D_j$ is 
$\tau$-quasi-uniform and $D_j$ pales in size to
$D$ resulting in a much favorable minimal separation radius and therefore a workable estimated value for  $\sigma_{\min}(\Phi_{D_j})$, which translates directly into superior algorithm execution stability. To further illustrate the idea behind DKI, we carried out several rounds of simulations\footnote{We employed a parallel computing technique in these simulations. More specifically, interpolation on $D_j$ for different $j$'s are executed by different PCs, which we refer to  as ``local machines".}; see the results 
in Figure \ref{Fig:dis-condition}. 

\begin{theorem}\label{Theorem:DKI-random}
Let  $1 \le m \le |D|.$ Suppose that $ \phi, \psi$ and $\varphi$ are SBDs satisfying the conditions set forth in Lemma \ref{Lemma:Narcowich-interpolation}. If $f^*\in N_\varphi$ and  $\overline{f}_D$ is the DKI estimator defined by \eqref{DKI} and based on the noisy data as in \eqref{Model1:fixed}, 
then
\begin{eqnarray}\label{dist-error}
         \mathbf E[\|\overline{f}_D-f^*\|_{\psi}^2] 
          &\leq &
         \frac{ C_3}{|D|} \sum_{j=1}^m 
         |D_j|^{\frac{1+(2\beta-2\alpha)\gamma}d}\|f^*\|^2_\varphi\\
         &+& \frac{ C_3M^2}{|D|^2} \sum_{j=1}^m     |D_j|^{\frac{4\gamma+2\gamma\beta+2+2d}d}   
           ,   \nonumber
\end{eqnarray}
where $C_3$ is a constant depending only  on $\alpha,\beta,\gamma,\tau$ and $d$.
\end{theorem}

We deem the two sums on the right hand sides of \eqref{dist-error} respectively as (deterministic) fitting-error estimate and (stochastic) sampling-error estimates. As opposed to KI, DKI gains strength in the sampling-error estimate while giving away some power in fitting-error estimate. In the grand scheme of things, this is a typical bias-variance trade-off. By choosing an optimal $m,$ DKI is capable of handling (statistically) data with large noise as the following result shows. 

\begin{definition}
   Let $0<\nu<1$ and $\tau>2$ be given. An input set $\Lambda$ is called a $(\nu, \tau)$-decomposable if  $\Lambda$  is the union of $m \sim |\Lambda|^\nu$ disjoint subsets $\Lambda_1, \ldots, \Lambda_m$ satisfying: (i) Each $\Lambda_j$ is $\tau$-quasi-uniform; (ii) $|\Lambda_j| \sim |\Lambda_{j'}|$ for all $j \ne j'.$
\end{definition}

\begin{corollary}\label{corollary:noise-data}
Let   $\hat{\phi}_k\sim k^{-2\gamma}$ with $2\gamma>d$. Suppose that the input set $\Lambda$ is $(\nu, \tau)$-decomposable with
$\nu={\frac{6\gamma+2}{6\gamma+d+2}}$. If $f^*\in\mathcal N_\phi$, then we have the following estimate 
$$
         \mathbf E[\|\overline{f}_D-f^*\|_{L^2(\mathbb S^d)}^2] 
          \leq 
         C_4 |D|^{-\frac{2\gamma}{6\gamma+d+2}}(M^2+\|f^*\|_\phi^2), \quad f^*\in\mathcal N_\phi, 
$$
where $C_4$ is a constant depending only on $\gamma,\tau$ and $d$.
\end{corollary}

In a sharp contrast to Proposition \ref{Corollary:inconsistence}, Corollary \ref{corollary:noise-data} shows the consistency of DKI in terms of handling noisy data of non-negligible magnitude and improves the approximation performance of KI in challenging computing environments.


\section{Operator Differences via  Quadrature Rules}\label{sec.sampling-ineq}

In this section, we first briefly exposit an integral-operator approach initiated in \cite{Feng2021radial} and then derive  tight upper bounds for  differences of operators of our interest, obtaining a certain type of Sobolev sampling inequalities \cite{hesse2021local} as a byproduct. Highlights of the section include Proposition \ref{Proposition:quadrature-for-convolution})  and Lemma \ref{Lemma:operator-difference}. 

\subsection{Integral operator approach and operator difference}
Given a positive definite kernel $\phi$, 
define $\mathcal L_\phi: L^2(\mathbb S^d)\rightarrow L^2(\mathbb S^d)$ to be the integral operator given by
$$
    \mathcal L_\phi f(x):=\int_{\mathbb{S}^{d}}\phi(x\cdot x')f(x')d\omega(x').
$$
The following lemma established in \cite{Feng2021radial} is a standard result to describe  relations among different SBDs via the integral operator $\mathcal L_\phi$.

\begin{lemma}\label{Lemma:Relations-among-Native}
Let  $\phi$,    $\psi$, $\varphi$ be SBDs satisfying \eqref{kernel-relation}.  Then, for any $f\in\mathcal N_\psi$  there holds
\begin{equation}\label{norm-relation}
    \|f\|_\psi =
         \|\mathcal L_\phi^{\frac{1-\beta}2}f\|_\phi. 
\end{equation}
Moreover, for $f^*\in\mathcal N_\varphi$, there exists an $h^*\in\mathcal N_\varphi$ such that
\begin{equation}\label{source-condition}
    f^*=\mathcal L_{\phi}^{\frac{\alpha-\beta}2}h^*,\qquad  \mbox{and} \quad
    \|f^*\|_\varphi=\|h^*\|_\psi,
\end{equation}
where $g(\mathcal L_\phi) f$ is defined by spectral calculus:
\[
     g(\mathcal L_\phi) f= \sum_{k=0}^\infty g(\hat{\phi}_k)\sum_{\ell=1}^{Z(d,k)}\hat{f}_{k,\ell} Y_{k,\ell}(x).
\]
\end{lemma}

Lemma \ref{Lemma:Relations-among-Native} shows the important role played by the integral operator $\mathcal L_\phi$ in terms of bridging native spaces associated with different SBDs. For the sake of clarity, we use $L_\phi$ to denote the restriction of  $\mathcal L_\phi$ on $\mathcal N_\phi$. 
Direct computation yields the following lemma.

\begin{lemma}\label{Lemma:integral-operator-relation}
Let $\kappa=\sup_{x\in\mathbb S^d}\sqrt{\phi(x\cdot x)}$. For any $g:[0,\kappa]\rightarrow\mathbb R_+$, there holds
\begin{equation}\label{integral-operator-relation}
    g(\mathcal L_\phi)f=g(L_\phi)f,\qquad \forall\ f\in\mathcal N_\phi.
\end{equation}
\end{lemma}

Let $\{w_i\}_{i=1}^{|D|}$ be a set of real numbers. Define $$S_{D,W}f=( \sqrt{w_1}f(x_1),\dots,\sqrt{w_{|D|}}f(x_{|D|}))^T$$
to be the weighted sampling operator from $\mathcal N_\phi$ to $\mathbb R^{|D|}$ \cite{smale2004shannon,smale2005shannon}. 
Write
$
    S^T_{D,W}{\bf c}=S_D^T{\sqrt{W}}{\bf c}= \sum_{i=1}^{|D|}\sqrt{w_{i}}c_i\phi_{x_i},
$
and 
\begin{equation}\label{weight-empirical-operaotr}
L_{\phi,D,W}f=S^T_{D,W}S_{D,W}f
=\sum_{i=1}^{|D|}w_{i}f(x_i)\phi_{x_i} 
\end{equation}
as the  positive operator from $\mathcal N_\phi$ to $\mathcal N_\phi$,
where $W={\rm diag} (w_{1},\dots,w_{|D|})$. 
As usual, we will refer to the above $w_i$'s as ``weights". If $w_1=\dots=w_{|D|}=1/|D|$, we denote $ S^T_{D,W}$ and   $L_{\phi,D,W}$ as $S_D^T$ and $L_{\phi,D}$ respectively. Since
$
    \frac1{|D|}\Phi_D=S_DS_D^T,
$
 it follows from \eqref{KI} that
\begin{equation}\label{operator-KI-1}
    f_D=S_D^T(S_DS_D^T)^{-1}  S_D {\bf y}_D=S_{D,W}^T(S_{D,W}S_{D,W}^T)^{-1}{\bf y}_{D,W},
\end{equation}
where ${\bf y}_{D,W}=(\sqrt{w_1}y_1,\dots,\sqrt{w_{|D|}}y_{|D|})^T.$

We call $\mathcal Q_{\Lambda,s}:=\{(w_{i,s},  x_i): w_{i,s}> 0
\hbox{~and~}   x_i\in \Lambda\}$
a positive  quadrature rule  of order $s$ on $\mathbb S^d$ if
\begin{equation}\label{eq:quadrature}
    \int_{\mathbb S^d}P(x)d\omega(x)=\sum_{x_i\in\Lambda} w_{i,s} P(  x_i), \qquad \forall P\in \mathcal P_{s}^d.
\end{equation}
 Positive spherical quadrature rules such as the ones in \eqref{eq:quadrature} are discussed in \cite{mhaskar2001spherical,brown2005approximation}. In particular, Theorem 3.1 \cite{brown2005approximation} asserts that for a $\tau$-quasi uniform point set $\Lambda$, there is a positive quadrature rule of order $s$ with $s\leq \bar{c}_1|\Lambda|^{1/d}$ and $0<w_{i,s}\leq \bar{c}_2|\Lambda|^{-1}$ for some constants $\bar{c}_1,\bar{c}_2$ depending only on $\tau,d,s$, which we will refer to as a D-type quadrature rule for short, as it plays a crucial role in our analysis in the present paper. For a D-type quadrature rule, we denote: 
\[
W_s:={\rm diag}(w_{1,s}, \dots, w_{|\Lambda|,s}) \quad {\rm and } \quad W^{1/2}_s:={\rm diag}(w^{1/2}_{1,s}, \dots, w^{1/2}_{|\Lambda|,s}).
\]  


 

  Our main result in this section is the following theorem, which gives   upper bounds for norms of products of pertinent operators when the weights in \eqref{weight-empirical-operaotr} are positive.

\begin{theorem}\label{Theorem:Product}
Let  $\mathcal Q_{\Lambda,s}:=\{(w_{i,s},  x_i): w_{i,s}> 0
\hbox{~and~}   x_i\in \Lambda\}$  be  a D-type
 quadrature rule of order $s$  on $\mathbb S^d$.
If $\hat \phi_k\sim k^{-2\gamma}$ with $\gamma> d/2$,  
then there is a constant $C'$  depending only on $\gamma$ and $d$, such that for any $\lambda\geq C' s^{-2\gamma}$, there holds
\begin{equation}\label{product-2}
     \|(L_{\phi,D,W_s}+\lambda
         I)^{-1/2}(L_\phi+\lambda I)^{1/2}\|\leq 2\sqrt{3}/3 
\end{equation}
and
\begin{equation}\label{product-3}
     \|(L_{\phi,D,W_s}+\lambda
         I)^{1/2}(L_\phi+\lambda I)^{-1/2}\|\leq \sqrt{6}/2, 
\end{equation}
 where $g(L_{\phi,D,W_s})$ and $g(L_\phi)$ are also defined by spectral calculus for $g: R_+ \rightarrow \mathbb R$.
\end{theorem}

The proof of Theorem \ref{Theorem:Product} is involved and will be given in the next subsection. Here we emphasize its usefulness in deriving approximation error estimates for KI's \eqref{minimal-norm-interpolation} and DKI's \eqref{DKI}. In particular, we use Theorem \ref{Theorem:Product} to prove a class of Sobolev-type sampling inequalities.


\begin{corollary}\label{corollary:sampling}
Let $\Lambda$ be a $\tau$-quasi-uniform set and $\mathcal Q_{\Lambda,s}:=\{(w_{i,s},  x_i): w_{i,s}> 0
\hbox{~and~}   x_i\in \Lambda\}$  a D-type
 quadrature rule  of order $s$ on $\mathbb S^d$. 
Suppose that $ \phi, \psi$ are SBDs satisfying $\hat{\phi} \sim k^{-2\gamma}$ with $\gamma > d/2$ and $\hat{\psi} \sim (\hat{\phi})^\beta$ for some $0 \le \beta \le 1.$
 Then there exists a constant  $C^*$  depending only on $\beta, \gamma$, and $d,$ such that for any 
  $\lambda\geq C^*  s^{-2\gamma}$ and 
 any $f\in\mathcal N_\phi$, there holds
$$
   \|f\|_{\psi}  
    \leq (4/3)^{(1-\beta)/2}\left(\lambda^{(1-\beta)/2}\|f\|_\phi+
    \lambda^{-\beta/2}\left(\sum_{i=1}^{|D|}w_{i,s}(f(x_i))^2\right)^{1/2}\right).
$$
\end{corollary}

For a quasi-uniform point set $\Lambda$, we have  $h_\Lambda\sim q_{\Lambda}\sim {|D|}^{-1/d}$.
Corollary \ref{corollary:sampling} then yields the following sampling inequality
$$
   \|f\|_\psi  
    \leq C_2^*\left(\lambda^{1/2}\|f\|_\phi+
    \left(\frac1{|D|}\sum_{i=1}^{|D|} (f(x_i))^2\right)^{1/2}\right),
$$
where $C_2^*$ is a constant independent of $|D|$. In particular, the case $\beta=0$ gives rise to the $L^2$ sampling inequality on spheres  established in  \cite[Theorem 5.1]{hesse2017radial}. We caution, however, that this implication is valid only under the assumption that the point set $\Lambda$ is quasi-uniform.

\begin{proof}[Proof of Corollary \ref{corollary:sampling}]
For a given $f\in \mathcal N_\phi$ and $0\leq \beta\leq 1$, we apply  Lemma \ref{Lemma:Relations-among-Native} and Theorem \ref{Theorem:Product} to get 
\begin{eqnarray*}
      \|f\|_{\psi} 
      &=&
      \|  L_\phi^{(1-\beta)/2}f\|_\phi \\
     & \leq &
     \| ( L_\phi+\lambda I)^{(1-\beta)/2}(L_{\phi,D,W_s}+\lambda I)^{-(1-\beta)/2}\|\|( L_{\phi,D,W_s}+\lambda I)^{(1-\beta)/2}f\|_\phi \\
      &\leq &
       (4/3)^{(1-\beta)/2}\lambda^{-\beta/2}\|( L_{\phi,D,W_s}+\lambda I)^{1/2}f\|_\phi\\
      &\leq &
         (4/3)^{(1-\beta)/2}\lambda^{-\beta/2}\| L_{\phi,D,W_s}^{1/2}f\|_\phi+   (4/3)^{(1-\beta)/2}\lambda^{-(1-\beta)/2}\|f\|_\phi.
\end{eqnarray*}
In developing the above inequalities, we have also used 
 the basic inequality $(a+b)^{1/2}\leq a^{1/2}+b^{1/2}$  for nonnegative real numbers $a,b$ and the well-known Cordes inequality \cite{bhatia2013matrix}  for positive operators $A,B$, 
\begin{equation}\label{Codes-inequality}
    \|A^u B^u\|\leq \|AB\|^u,\qquad 0\leq u\leq 1. 
\end{equation}
By the definition of $L_{\phi,D,W_s}$, we deduce
\begin{eqnarray*}
  &&\|L_{\phi,D,W_s}^{1/2}f\|_\phi^2
  =
  \langle L_{\phi,D,W_s}^{1/2}f,L_{\phi,D,W_s}^{1/2}f\rangle_\phi
  =
  \langle f,L_{\phi,D,W_s}f\rangle_\phi 
  = 
  \langle f,S_D^TW_sS_Df\rangle_\phi\\
  &=&
  \langle  \sqrt{W_s}S_D f,\sqrt{W_s}S_Df\rangle_\phi
  =\sum_{i=1}^{|D|}w_{i,s}(f(x_i))^2.
\end{eqnarray*}
It then follows that
\begin{eqnarray*}
 \|f\|_{\psi}  
   &\leq& 
   (4/3)^{(1-\beta)/2}\lambda^{-\beta/2}\left(\sum_{i=1}^{|D|} w_{i,s}(f(x_i))^2\right)^{1/2}+   (4/3)^{(1-\beta)/2}\lambda^{-(1-\beta)/2}\|f\|_\phi,
\end{eqnarray*}
which is the desired result. 
\end{proof}

\subsection{Quadrature rule for products of functions}

To prove Theorem \ref{Theorem:Product}, we need 
a quadrature rule for products of two functions in $\mathcal N_\phi$ associated with certain positive operators of our interest.

\begin{proposition}\label{Proposition:quadrature-for-convolution}
Let $\hat{\phi}_k\sim k^{-2\gamma}$ with $\gamma>d/2$ and $u,v\in[0,1/2]$. If $\Lambda=\{  x_i\}_{i=1}^{|\Lambda|}$ is $\tau$-quasi uniform and $s\leq \tilde{c}_1|\Lambda|^{1/d}$, then  for a D-type  quadrature  rule $\mathcal Q_{\Lambda,s}=\{(w_{i,s},  x_i):
\hbox{~and~}   x_i\in \Lambda\}$  of order $s$ and
each pair of $f, g \in \mathcal N_\phi$, there holds
\begin{align*} 
   &\left|\int_{\mathbb{S}^{d}}  (\eta_{\lambda,u}( L_\phi)f)  (x)(\eta_{\lambda,v}( L_\phi)g)  (x) d \omega(x)-\sum_{x_i\in\Lambda} w_{i, s} \left[(\eta_{\lambda,u}( L_\phi)f)  (x_{i} )(\eta_{\lambda,v}( L_\phi)g) (x_i )\right]\right|\\
   &\leq 
  \tilde{c}'\lambda^{-\max\{u,v\}}s^{-\gamma}
   \|f\|_\phi\|g\|_\phi,
\end{align*}
 in which   $\eta_{\lambda,u}( L_\phi)f=(L_\phi+\lambda I)^{-u}f$ for $0\leq u\leq 1/2$  and $0<\lambda\leq 1$,  and  $\tilde{c}' $  is a constant depending only on $\gamma,\kappa,u,v,\tau$ and $d$.  
\end{proposition}

Proposition \ref{Proposition:quadrature-for-convolution} extends  \cite[Corollary 3.6]{brauchart2007numerical} which asserts that 
\begin{equation}\label{fixed-cubabuture-native}
    \left|\int_{\mathbb S^d}f(x)d\omega(x)-\sum_{x_i\in\Lambda} w_{i,s} f(  x_{i})\right|\leq \tilde{c}_1's^{-\gamma}\|f\|_{\phi}, \quad f \in \mathcal N_\phi,
\end{equation}
where $\tilde{c}_1'$ is a constant depending only on $\gamma,d$.
 As a matter of fact, setting $u=v=0$ and $g(x)\equiv1$ in Proposition \ref{Proposition:quadrature-for-convolution} recovers \eqref{fixed-cubabuture-native}. Furthermore, the results of Proposition \ref{Proposition:quadrature-for-convolution}  improve  those of \cite[Proposition 1]{Feng2021radial}, which gives the following estimate:
 \begin{align*}
   &\left|\int_{\mathbb{S}^{d}}  (\eta_{\lambda,u}( L_\phi)f)  (x)(\eta_{\lambda,v}( L_\phi)g)  (x) d \omega(x)-\sum_{x_i\in\Lambda} w_{i, s} \left[(\eta_{\lambda,u}( L_\phi)f)  (x_{i} )(\eta_{\lambda,v}( L_\phi)g) (x_i )\right]\right| \\
  & \leq  \tilde{c}\|f\|_\phi\|g\|_\phi
   \begin{cases}
   \lambda^{-u-v}s^{-2\gamma}+s^{-(1-u-v)\gamma},    &\mbox{if}\ u+v<1, \\
     \lambda^{-u-v}s^{-2\gamma}+\log s,    &\mbox{if}\ u+v=1.
      \end{cases} 
      \end{align*}

Obtaining the sharper upper bound here requires the following two lemmas from respectively Narcowich et al \cite[Corollary 5.4]{narcowich2007direct} and Dai \cite{dai2006multivariate}.



\begin{lemma}\label{Lemma:polynomial-interpolation}
Let $\hat{\phi}_k\sim k^{-2\gamma}$ with $\gamma>d/2$. Let $\Lambda$ be a finite subset of $\mathbb S^d$. Then there exists a constant $\tilde{c}_1$ depending only on $\gamma$ and $d$, such that for any $f \in\mathcal N_\phi$ and $s^*\geq \tilde{c}_1/q_\Lambda$, there is a $P_{s^*}\in\mathcal P_{s^*}^d$ satisfying $\|P_{s^*}\|_\phi\leq 6\|f\|_\phi$ and $\|P_{s^*}\|_{L^2(\mathbb S^d)}\leq 6\|f\|_{L^2(\mathbb S^d)}$.
Moreover, we have
\begin{equation}\label{best-app-interpolation1}
     f(x_i)=P_{s^*}(x_i), \qquad i=1,\dots,|\Lambda|,
\end{equation}
and
\begin{equation}\label{best-app-interpolation2}
    \|f-P_{s^*}\|_{L^2(\mathbb S^d)}\leq \tilde{c}_2h_\Lambda^{\gamma} \|f\|_\phi,
\end{equation}
where $\tilde{c}_2$ is a constant depending only on $d,\beta$ and $\gamma$.
\end{lemma}


\begin{lemma}\label{Lemma:MZ-inequality}
Let $s\in\mathbb N$ and $\Lambda$ be quasi-uniform. Suppose that $\mathcal Q_{\Lambda,s}:=\{(w_{i,s},  x_i): x_i\in \Lambda\}$ is a  D-type  quadrature rule of order $s$  on $\mathbb S^d$. Then for   any  $P\in \Pi_{s'}^d$ with $s'\in\mathbb N$, there holds
 \begin{eqnarray*}
   \sum_{x_i\in\Lambda}w_{i,s} |P(x_i)|^2  &\leq& \tilde{c}_1(s'/s)^{d} \|P\|^2_{L^2(\mathbb S^d)},\qquad s'>\lfloor s/3\rfloor,\\
\end{eqnarray*}
where $\tilde{c}_1$ is a constant  depending only on $d$.
\end{lemma}

 The above lemma can be 
 derived by combining \cite[Lemma 4.4]{dai2006multivariate} and \cite[Lemma 4.5]{dai2006multivariate}  (see also \cite[Theorem 2.1]{dai2006generalized} and  \cite[Theorem 3.3]{mhaskar2006weighted}). 
With the help of above lemmas, we are in a position to prove Proposition \ref{Proposition:quadrature-for-convolution}.

\begin{proof}[Proof of Proposition \ref{Proposition:quadrature-for-convolution}]
Let $\tilde{c}\geq \tilde{c_1}$ be a positive constant properly chosen so that $s:=\tilde{c}/q_\Lambda$ satisfies the pertinent interpolation condition in Lemma \ref{Lemma:polynomial-interpolation}. Let $P_f^*$ and $P_g^*$ be, respectively, the  interpolating polynomials (in accordance to Lemma \ref{Lemma:polynomial-interpolation}) from $\Pi^d_{s}$ 
 for $\eta_{\lambda,u}( L_\phi)f$ and $\eta_{\lambda,v}( L_\phi)g\in\mathcal N_\phi$. 
Using a standard ``zero-sum" technique and  H\"{o}lder's inequality, we derive 
that
\begin{align*} 
    &\left|\int_{\mathbb{S}^{d}}  (\eta_{\lambda,u}( L_\phi)f)  (x)(\eta_{\lambda,v}( L_\phi)g)  (x) d \omega(x)-\sum_{x_i\in\Lambda} w_{i, s} \left[(\eta_{\lambda,u}( L_\phi)f)  (x_{i} )(\eta_{\lambda,v}( L_\phi)g) (x_i )\right]\right|\\
    &\leq
  \mathcal A_1+\mathcal A_2+\mathcal A_3,  
\end{align*}
in which
\begin{eqnarray*}
  \mathcal A_1&:=&\| (\eta_{\lambda,u}( L_\phi)f)\|_{L^2(\mathbb S^d)}\|(\eta_{\lambda,v}( L_\phi)g)-P_g^*\|_{L^2(\mathbb S^d)},\\
  \mathcal A_2&:=&
  \|P_g^*\|_{L^2(\mathbb S^d)}\|(\eta_{\lambda,u}( L_\phi)f)-P_f^*\|_{L^2(\mathbb S^d)},\\
  \mathcal A_3&:=&
  \left|\int_{\mathbb{S}^{d}}  P_f^*(x)P_g^*(x) d \omega(x)-\sum_{x_i\in\Lambda} w_{i, s}P_f^*(x_i)P_g^*(x_i)\right|.
\end{eqnarray*}
 Since $0\leq u,v\leq 1/2$, we get from \eqref{norm-relation} that 
$$
  \| (\eta_{\lambda,u}( L_\phi)f)\|_{L^2(\mathbb S^d)}\|\leq \kappa^{1-2u}\|f\|_\phi,\quad    \| (\eta_{\lambda,v}( L_\phi)g)\|_{L^2(\mathbb S^d)}\|\leq \kappa^{1-2v}\|g\|_\phi.
$$
We use the result of Lemma \ref{Lemma:polynomial-interpolation} to bound $\mathcal A_1$ and $\mathcal A_2$ as follows.
 \begin{equation}\label{bound-A1}
   \mathcal A_1\leq \tilde{c}_2\kappa^{1-2u}h_\Lambda^{\gamma} \|\eta_{\lambda,v}( L_\phi)g\|_\phi\|f\|_\phi\leq \tilde{c}_2\kappa^{1-2u}\lambda^{-v}h_\Lambda^{\gamma} \|f\|_\phi\|  g\|_\phi, \quad \mbox{and}
\end{equation}
\begin{equation}\label{bound-A2}
   \mathcal A_2\leq 6\tilde{c}_2\kappa^{1-2v}h_\Lambda^{\gamma} \|\eta_{\lambda,u}( L_\phi)f\|_\phi\|g\|_\phi\leq 6\tilde{c}_2\kappa^{1-2v}\lambda^{-u}h_\Lambda^{\gamma} \|f\|_\phi\|  g\|_\phi.
\end{equation}
It remains to bound $\mathcal A_3$. To this end,  we let
 $Q^*_1,Q_2^*\in\mathcal P_{\lfloor s/6\rfloor}^d$ be the $L^2$-projections of $P_f^*$ and $P_g^*$ onto 
 $\mathcal P_{\lfloor s/6\rfloor}^d$. The the following inequalities are immediate.
\begin{align}\label{best-app-for-p-rate}
   &\|P_f^*-Q_1\|_{L^2(\mathbb S^d)}\leq \tilde{c}_3s^{-\gamma}\|\eta_{\lambda,u}( L_\phi)f\|_\phi,\quad {\rm and} \\ \nonumber
   &\|P_g^*-Q_2\|_{L^2(\mathbb S^d)}\leq \tilde{c}_3s^{-\gamma}\|\eta_{\lambda,v}( L_\phi)g\|_\phi,
\end{align}
where $\tilde{c}_3$ is a constant depending only on $\gamma$. We then break up $\mathcal A_3$ into five terms and use H\"{o}lder's inequality to write
$
   \mathcal A_3\leq \sum_{k=1}^5\mathcal A_{3,k},
$
where
\begin{eqnarray*}
  \mathcal A_{3,1}&:=& \|P^*_f-Q_1\|_{L^2(\mathbb S^d)}\|P^*_g\|_{L^2(\mathbb S^d)},  \\
  \mathcal A_{3,2} &:=&   \|P^*_g-Q_2\|_{L^2(\mathbb S^d)}\|Q_1\|_{L^2(\mathbb S^d)},   \\
  \mathcal A_{3,3} &:=& \left(\sum_{x_i\in\Lambda}w_{i,s}(P^*_f(x_i)-Q_1(x_i))^2\right)^{1/2}\left(\sum_{x_i\in\Lambda}w_{i,s}(P^*_g(x_i))^2\right)^{1/2},\\
  \mathcal A_{3,4} &:=& \left(\sum_{x_i\in\Lambda w_{i,s}}w_{i,s}(P^*_g(x_i)-Q_2(x_i))^2\right)^{1/2}\left(\sum_{x_i\in\Lambda}w_{i,s}(Q_1(x_i))^2\right)^{1/2},\\
   \mathcal A_{3,5} &:=&  \left|\int_{\mathbb{S}^{d}}  Q_1(x)Q_2(x) d \omega(x)-\sum_{x_i\in\Lambda} w_{i, s}Q_1(x_i)Q_2(x_i)\right|.
\end{eqnarray*}
The fact that $Q_1Q_2\in\mathcal P_s^d$ implies  $\mathcal A_{3,5}=0$. Furthermore, 
 \eqref{best-app-for-p-rate} together with 
$\|Q_1\|_{L^2(\mathbb S^d)}\leq \|P_f^*\|_{L^2(\mathbb S^d)}$ and $\|Q_2\|_{L^2(\mathbb S^d)}\leq \|P_g^*\|_{L^2(\mathbb S^d)}$ yields
\begin{align*}
&   \mathcal A_{3,1}\leq  6\tilde{c}_3s^{-\gamma}\|\eta_{\lambda,u}( L_\phi)f\|_\phi\|\eta_{\lambda,v}(L_\phi)g\|_{L^2(\mathbb S^d)}\leq 6\tilde{c}_3\kappa^{1-2v}\lambda^{-u}s^{-\gamma}
   \|f\|_\phi\|g\|_\phi, \; {\rm and}\\
   &\mathcal A_{3,2}\leq  6\tilde{c}_3\kappa^{1-2u}\lambda^{-v}s^{-\gamma}\|f\|_\phi\|g\|_\phi.  
\end{align*}
Using Lemma \ref{Lemma:MZ-inequality} with $s'=s^*=\tilde{c}q_\Lambda^{-1}$, we have 
\begin{align*}
&   A_{3,3}\leq \tilde{c}_1(\tilde{c}q_\Lambda^{-1}s^{-1})^{ 2d} \mathcal A_{3,1}\leq 6\tilde{c}_1(\tilde{c}q_\Lambda^{-1}s^{-1})^{2d}\tilde{c}_3\kappa^{1-2v}\lambda^{-u}s^{-\gamma}
   \|f\|_\phi\|g\|_\phi, \; {\rm and}\\
&   A_{3,4}\leq \tilde{c}_1(\tilde{c}q_\Lambda^{-1}s^{-1})^{ 2d} \mathcal A_{3,2}\leq 6\tilde{c}_1(\tilde{c}q_\Lambda^{-1}s^{-1})^{2d}\tilde{c}_3\kappa^{1-2u}\lambda^{-v}s^{-\gamma}
   \|f\|_\phi\|g\|_\phi.  
\end{align*}
Combining  the above five estimates, we have
\begin{equation}\label{Bound-A3}
   \mathcal A_{3}\leq
    12(\tilde{c}_1(\tilde{c}q_\Lambda^{-1}s^{-1})^{2d}\tilde+2)
    \tilde{c}_3(\kappa+1)\lambda^{-\max\{u,v\}}s^{-\gamma}
   \|f\|_\phi\|g\|_\phi. 
\end{equation}
We then get from \eqref{bound-A1}, \eqref{bound-A2} and \eqref{Bound-A3}  that
\begin{align*}
  &\left|\int_{\mathbb{S}^{d}}  (\eta_{\lambda,u}( L_\phi)f)  (x)(\eta_{\lambda,v}( L_\phi)g)  (x) d \omega(x)-\sum_{x_i\in\Lambda} w_{i, s} \left[(\eta_{\lambda,u}( L_\phi)f)  (x_{i} )(\eta_{\lambda,v}( L_\phi)g) (x_i )\right]\right|\\
 & \leq
  12(\tilde{c}_1(\tilde{c}q_\Lambda^{-1}s^{-1})^{2d}\tilde+2)
    \tilde{c}_3(\kappa+1)\lambda^{-\max\{u,v\}}s^{-\gamma}
   \|f\|_\phi\|g\|_\phi\\
  &\;\;+ 7\tilde{c}_2(\kappa+1)\lambda^{-\max\{u,v\}}h_\Lambda^{\gamma} \|f\|_\phi\|  g\|_\phi.
\end{align*}
 Since
  $\Lambda$ is $\tau$-quasi-uniform, we have 
 $
  q_\Lambda\sim |D|^{-1/d}\sim h_\Lambda\sim s^{-1}.
 $
It then follows that
\begin{align*}
  &\left|\int_{\mathbb{S}^{d}}  (\eta_{\lambda,u}( L_\phi)f)  (x)(\eta_{\lambda,v}( L_\phi)g)  (x) d \omega(x)-\sum_{x_i\in\Lambda} w_{i, s} \left[(\eta_{\lambda,u}( L_\phi)f)  (x_{i} )(\eta_{\lambda,v}( L_\phi)g) (x_i )\right]\right|\\
  &\leq 
  \tilde{c}'\lambda^{-\max\{u,v\}}s^{-\gamma}
   \|f\|_\phi\|g\|_\phi,
\end{align*}
where $\tilde{c}'$ is a constant depending only on $\tau,u,v,\kappa$, $\gamma$ and  $d$.
This completes the proof of Proposition \ref{Proposition:quadrature-for-convolution}.
\end{proof}

\subsection{Operator differences and products} 

Here we first use Proposition \ref{Proposition:quadrature-for-convolution} to bound $\|L_\phi - L_{\phi,D,W_s}\|$ and  then prove Theorem \ref{Theorem:Product}. 



\begin{lemma}\label{Lemma:operator-difference}
Let $\mathcal Q_{\Lambda,s}:=\{(w_{i,s},  x_i): w_{i,s}>0
\hbox{~and~}   x_i\in \Lambda\}$ be  a positive
 quadrature rule   on $\mathbb S^d$ with degree $s\in\mathbb N$.
If $\hat \phi_k\sim k^{-2\gamma}$  with $\gamma> d/2$,
then for any  $\lambda \geq s^{-2\gamma}$ and $u,v\in[0,1/2]$, there holds
\begin{eqnarray}\label{operaotr-difference-111}
  \|(  L_\phi+\lambda I )^{-u}(L_{\phi,D,W_s}-  L_\phi)(  L_\phi+\lambda I)^{-v}\|  
  \leq  \tilde{c}'\lambda^{-\max\{u,v\}}s^{-\gamma}.
\end{eqnarray}
\end{lemma}

\begin{proof}
Due to the definition  of operator norm, we have from the reproducing property of $\phi$, i.e. $\langle\phi_x,f\rangle_\phi=f(x)$ for any $f\in\mathcal N_\phi$  that
\begin{eqnarray*}
  &&\|(  L_\phi+\lambda I)^{-u}(L_{\phi,D,W_s}-  L_\phi)(  L_\phi+\lambda I)^{-v}\|\\
  &=&
  \sup_{\|f\|_\phi\leq 1}\|(  L_\phi+\lambda I)^{-u}(L_{\phi,D,W_s}-  L_\phi)(  L_\phi+\lambda I)^{-v}f\|_\phi\\
  &=&
  \sup_{\|f\|_\phi\leq 1}\sup_{\|g\|_\phi\leq 1}
  \langle (  L_\phi+\lambda I)^{-u}(L_{\phi,D,W_s}-  L_\phi)(  L_\phi+\lambda I)^{-v}f,g\rangle_\phi\\
  &=&
  \sup_{\|f\|_\phi\leq 1}\sup_{\|g\|_\phi\leq 1}\langle (L_{\phi,D,W_s}- L_\phi)(  L_\phi+\lambda I)^{-u}f,
  (  L_\phi+\lambda I)^{-v}g\rangle_\phi\\
  &=&
  \sup_{\|g\|_\phi\leq 1,\|f\|_\phi\leq 1}
  \left|\left\langle\int_{\mathbb S^d}(  L_\phi+\lambda I)^{-u}f(x')\phi_{x'}d\omega(x')
  \right.\right.\\
  &-&
  \left.\left.\sum_{i=1}^{|D|}w_{i,s}(  L_\phi+\lambda I)^{-u}f(x_i)\phi_{x_i}, (  L_\phi+\lambda  I )^{-v}g\right\rangle_\phi\right|\\
  &=&
  \sup_{\|g\|_\phi\leq 1,\|f\|_\phi\leq 1} \left|\int_{\mathbb S^d}(  L_\phi+\lambda I)^{-u}f(x')\langle\phi_{x'},(  L_\phi+\lambda I)^{-v} g\rangle_\phi d\omega(x')\right.\\
  &-& \left. \sum_{i=1}^{|D|}w_{i,s}(  L_\phi+\lambda I)^{-u}f(x_i)\langle\phi_{x_i}, (  L_\phi+\lambda I)^{-v}g\rangle_\phi\right|\\
  &=&
  \sup_{\|g\|_\phi\leq 1,\|f\|_\phi\leq 1}\left|\int_{\mathbb S^d}(  L_\phi+\lambda I)^{-u}f(x')(  L_\phi+\lambda I )^{-v}g(x')d\omega(x')\right.\\
  &-&
  \left.\sum_{i=1}^{|D|}w_{i,s}(  L_\phi+\lambda I)^{-u}f(x_i)(  L_\phi+\lambda I)^{-v}g(x_i)\right|.
\end{eqnarray*}
Then, \cref{Proposition:quadrature-for-convolution}  implies
$$
     \|(  L_\phi+\lambda I)^{-1/2}(L_{\phi,D,W_s}-  L_\phi)(  L_\phi+\lambda I)^{-1/2}\|\leq \tilde{c}'\lambda^{-\max\{u,v\}}s^{-\gamma}. 
$$
This completes the proof of   \cref{Lemma:operator-difference}.
\end{proof}


\begin{proof}[Proof of Theorem \ref{Theorem:Product}]
For  $\lambda\geq 16\tilde{c}'^2 s^{-2}$, Lemma \ref{Lemma:operator-difference} yields
\begin{equation}\label{R11111}
   \|(  L_\phi+\lambda I )^{-1/2}(L_{\phi,D,W_s}-  L_\phi)(  L_\phi+\lambda I)^{-1/2}\| \leq 1/4.     
\end{equation}
 Then,
\begin{align*}
      &(L_\phi+\lambda I)^{1/2}(L_{\phi,D,W_s}+\lambda
         I)^{-1}(L_\phi+\lambda I)^{1/2}\\
         =&I +
         (L_\phi+\lambda I)^{1/2}[(L_{\phi,D,W_s}+\lambda
         I)^{-1}-(L_{\phi}+\lambda
         I)^{-1}]   (L_\phi+\lambda I)^{1/2} \\
          =& I+ 
         (L_\phi+\lambda I)^{-1/2}(L_\phi-L_{\phi,D,W_s})(L_\phi+\lambda I)^{-1/2}(L_\phi+\lambda I)^{1/2} \\
         &\times (L_{\phi,D,W_s}+\lambda
         I)^{-1}(L_\phi+\lambda I)^{1/2}.
\end{align*}
Thus,
\begin{eqnarray*}
  && \|(L_\phi+\lambda I)^{1/2}(L_{\phi,D,W_s}+\lambda
         I)^{-1}(L_\phi+\lambda I)^{1/2}\|\\
         &\leq&
         1+ \frac14\|(L_\phi+\lambda I)^{1/2}(L_{\phi,D,W_s}+\lambda
         I)^{-1}(L_\phi+\lambda I)^{1/2}\|.
\end{eqnarray*}
This implies
$$
  \|(L_\phi+\lambda I)^{1/2}(L_{\phi,D,W_s}+\lambda
         I)^{-1}(L_\phi+\lambda I)^{1/2}\|\leq \frac43,
$$
and therefore proves \eqref{product-2}. We then turn to proving \eqref{product-3}.  It follows from \eqref{product-2} and \eqref{R11111} that 
\begin{align*}
    &\|(L_{\phi,D,W_s}+\lambda I)^{-1/2}  (L_{\phi,D,W_s}- L_{\phi})(L_{\phi,D,W_s}+\lambda I)^{-1/2}\|\\
    \leq&
    \|(L_{\phi,D,W_s}+\lambda I)^{-1/2}(L_{\phi}+\lambda I)^{1/2}\|^2 \\ 
    & \times \|(L_{\phi}+\lambda I)^{-1/2}(L_{\phi,D,W_s}- L_{\phi})(L_{\phi}+\lambda I)^{-1/2}\|
    \leq 
    \frac13.
\end{align*}
Therefore,
\begin{eqnarray*}
    &&\|(L_{\phi,D,W_s}+\lambda I)^{1/2}(L_{\phi}+\lambda
         I)^{-1} (L_{\phi,D,W_s}+\lambda I)^{1/2}\|\\ 
        &\leq &
   1+ \frac13\|(L_{\phi,D,W_s}+\lambda I)^{1/2}(L_{\phi}+\lambda
         I)^{-1} (L_{\phi,D,W_s}+\lambda I)^{1/2}\|,
\end{eqnarray*}
which implies 
$$
   \|(L_{\phi,D,W_s}+\lambda I)^{1/2}(L_{\phi}+\lambda
         I)^{-1} (L_{\phi,D,W_s}+\lambda I)^{1/2}\|\leq \frac32.
$$
This proves \eqref{product-3} and
completes the proof of Theorem \ref{Theorem:Product}.
\end{proof}

\section{Proofs}\label{Sec.Proofs}
This section provides proofs of our main results.

\begin{proof}[Proof of Proposition \ref{Corollary:inconsistence}]
In view of \eqref{KI}, we use the reproducing property of $\phi$ to write
$$
   \|f_D\|_\phi^2=\left\langle \sum_{i=1}^{|D|}a_i\phi_{x_i},\sum_{j=1}^{|D|}a_j\phi_{x_j}\right\rangle_\phi
   =\sum_{i=1}^{|D|}\sum_{j=1}^{|D|}a_ia_j\phi(x_i\cdot x_j)={\bf a}_D^T\Phi_D{\bf a}_D,
$$
in which ${\bf a}_D=\Phi_D^{-1}{\bf y}_D$. Let $\sigma_1, \ldots, \sigma_{|D|}$ be the eigenvalues of $\Phi_D$ in descending order. The invariance of trace yields the identity: $\sigma_1 + \cdots + \sigma_{|D|} = |D| \phi(1).$
Thus we have that $\sigma_1  < |D| \phi(1),$ and therefore that $\xi \ge (|D| \phi(1))^{-1},$ where $\xi$ denotes the smallest eigenvalue of $\Phi_D^{-1}.$ Hence, 
$$
   \|f_D\|_\phi^2=({\bf y}_D)^T \Phi_D^{-1}{\bf y}_D \ge \xi \sum^{|D|}_{i=1} y_i^2 \ge (\phi(1))^{-1} \theta^2 M^2.
$$
This completes the proof.
\end{proof}

To prove Theorem \ref{Theorem:DKI-random},
we adopt  a standard error decomposition strategy in learning theory \cite{smale2005shannon,lin2017distributed} and the sampling inequality developed in Section \ref{sec.sampling-ineq}. Let  $G_D$ be the projection operator from $\mathcal N_\phi$ to
$$
   \mathcal H_{\phi,D}=\mbox{span}\{\phi_{x_1},\dots,\phi_{x_{|D|}}\}. 
$$ 
Then it is easy to get 
\begin{equation}\label{Projection}
     G_D:=S_D^T(S_D S_D^T)^{-1}S_D.
\end{equation}  
Define
\begin{equation}\label{def.noise-free}
     f_D^\diamond= G_Df^*
\end{equation}
be the noise-free version of $f_D$. Then the triangle inequality yields
\begin{equation}\label{first-error-dec}
    \|f_D-f^*\|_\psi\leq
    \|f_D^\diamond-f^*\|_\psi+
    \|f_D^\diamond-f_D\|_\psi.
\end{equation}
The above two terms are called as the approximation error and sample error, respectively.  To bound the approximation error, we need the following lemma whose proof is standard and can be eaisly derived from \cite{rudi2015less,lin2020kernel}.

\begin{lemma}\label{lemma:operator-bound-PD}
Let $u,v\geq 0$ and $\lambda>0$. If $\mathcal Q_{\Lambda,s}:=\{(w_{i,s},  x_i): w_{i,s}> 0
\hbox{~and~}   x_{i}\in \Lambda\}$ is 
a positive  quadrature rule   on $\mathbb S^d$ with degree $s\in\mathbb N$, then 
\begin{equation}\label{projection-bound}
      \|  L_\phi^{v}(I-G_D)L_\phi^u\|\leq \lambda \|  L_\phi^u(L_{\phi,D,W_s}+\lambda I)^{-1}L_\phi^v\|.
\end{equation}
\end{lemma}

The above lemma  describes the difference between $G_D$ and the identity mapping $I$. We remove the proof of Lemma  \ref{lemma:operator-bound-PD} for the sake of brevity. 
We then can bound the approximation error as follows.
 
\begin{lemma}\label{Lemma:Approx-error}
 If $f\in N_\varphi$ with $\varphi$ satisfying  (\ref{kernel-relation})  for $\alpha\geq 1$ and $0\leq \alpha-\beta\leq 1$ and  $\hat{\phi}_k\sim k^{-2\gamma}$ with $\gamma>d/2$, then  for some $s\sim |D|^{-1/d}$ and any $\lambda\geq 16 \tilde{c}'^2 s^{-2\gamma}$, there holds
 \begin{equation}\label{approx.error-1}
 \|f_D^\diamond-f^*\|_\psi
    \leq \frac43\lambda^{\frac{\alpha-\beta}2}\|f^*\|_\varphi
\end{equation}
where $\bar{c}_1$ is a constant depending only on $\tilde{c}$, $c$, $\alpha$,    $\beta$ and $\gamma$.
\end{lemma}

\begin{proof}
Since $f^*\in\mathcal N_\varphi$, we have from \eqref{source-condition} and Lemma \ref{Lemma:integral-operator-relation} that
\begin{eqnarray*}
  && \|f_D^\diamond-f^*\|_\psi
  =
  \|L_\phi^{\frac{1-\beta}2}(G_D-I)\mathcal L_{\phi}^{\frac{\alpha-\beta}2}h^*\|_\phi\\
  &\leq&
  \|L_{\phi}^{\frac{1-\beta}2}(G_D-I)L_{\phi}^{\frac{\alpha-1}2}\|  \|L_\phi^{\frac{1-\beta}2}h^*\|_\phi
  =
  \|L_{\phi}^{\frac{1-\beta}2}(G_D-I)L_{\phi}^{\frac{\alpha-1}2}\|\|h^*\|_\psi.
\end{eqnarray*}
But Theorem \ref{Theorem:Product} and  Lemma \ref{lemma:operator-bound-PD}   with $u=\frac{1-\beta}2$ and $v=\frac{1-\beta}2$
     imply that for any $\lambda\geq   16\tilde{c}'^2s^{-2\gamma}$   and $\alpha\leq \beta+1\leq 2$, there holds
\begin{eqnarray*}
   &&
   \|L_{\phi}^{\frac{1-\beta}2}(G_D-I)L_{\phi}^{\frac{\alpha-1}2}\|
   \leq 
    \lambda \|L_{\phi}^{\frac{1-\beta}2}(L_{\phi,D,W_s}+\lambda I)^{-1}L_{\phi}^{\frac{\alpha-1}2}\|\\
     &\leq&
     \lambda \|(L_{\phi}+\lambda I)^{\frac{1-\beta}2}(L_{\phi,D,W_s}+\lambda I)^{-1}(L_{\phi}+\lambda I )^{\frac{\alpha-1}2}\|\\
     &\leq&
     \frac43\lambda^{\frac{\alpha-\beta}2}
\end{eqnarray*}
Therefore,  we obtain \eqref{approx.error-1} by noting \eqref{source-condition} and then complete  the proof of Lemma \ref{Lemma:Approx-error}.
\end{proof}

To prove Theorem \ref{Theorem:DKI-random}, we  need the following lemma from \cite{Feng2021radial}.

\begin{lemma}\label{lemma:value-difference-random}
Let $\delta\in(0,1)$. If \eqref{Model1:fixed} holds, $\hat{\phi}_k\sim k^{-2\gamma}$ with $\gamma>d/2$ and $\{\varepsilon_i\}_{i=1}^{|D|}$ are a set of i.i.d. random variables satisfying $E[\varepsilon_i]=0$ and $|\varepsilon_i|\leq M$ for all $i=1,2,\dots,|D|$, then for any diagnosis matrix $W$ with  diagonal element being $\{w_{1},\dots,w_{|D|}\}$ and $0<w_i \leq c_1|D|^{-1} $,   with confidence $1-\delta$, there holds 
\begin{eqnarray}\label{norm-difference-random}
     \left\|(  L_\phi+\lambda I)^{-1/2} (L_{\phi,D,W}f^*-S_{D,W}^T{\bf y}_{D,W})\right\|_\phi
     \leq 
    \bar{c}_1M\lambda^{-\frac{d}{4\gamma}} |D|^{-1/2}\log\frac3\delta,
\end{eqnarray}
where $ \bar{c}_1$ is a constant  depending only on $d$, $\gamma$ and  $\|f^*\|_{\varphi}$.
\end{lemma}

In view of \eqref{first-error-dec}, it is sufficient to bound the approximation and sample errors, respectively. The approximation error has been handled in Lemma \ref{Lemma:Approx-error}. In what follows, we work on bounding the sample error $\|f_D^\diamond-f_D\|_\psi$.

\begin{lemma}\label{Lemma:Sample-error}
Let $\delta\in(0,1)$, $\hat{\phi}_k\sim k^{-2\gamma}$ and $\gamma>d/2$. If \eqref{Model1:fixed} holds with $\{\varepsilon_i\}_{i=1}^{|D|}$   a set of i.i.d. random variables satisfying $E[\varepsilon_i]=0$ and $|\varepsilon_i|\leq M$ for all $i=1,2,\dots,|D|$,   $f\in\mathcal N_\varphi$ with  (\ref{kernel-relation}), then for any  $ \mu\geq \bar{c}_2|D|^{-2\gamma/d}$, with confidence $1-\delta$, there holds
 \begin{equation}\label{sample.error-1-fixed}
    \|f_D^\diamond-f_D\|_\psi\leq
  \bar{c}_3 \mu^{-\frac\beta2}  (1+\mu |D| (\sigma_{|D|} (\Phi_D))^{-1}) 
  M\mu^{-\frac{d}{4\gamma}} |D|^{-1/2}\log\frac3\delta,  
\end{equation} 
where $\bar{c}_2,\bar{c}_3$ are constants depending only on $\tilde{c}'$, $d$ and $\gamma$.
\end{lemma}

\begin{proof}
For an arbitrary $\eta>0$, it follows from
\begin{equation}\label{elementary-compu-inverse-1}
     \frac1t=\frac1{t+\eta}+\frac{\eta}{(t+\eta)^2}+\frac{\eta^2}{t(t+\eta)^2},\qquad\forall t>0,\quad \eta>0,  
\end{equation}
that for any $\mu>0$ there hold
\begin{eqnarray}\label{exchang-inverse-matrix}
	(S_{D,W_s}S_{D,W_s}^T)^{-1}
		&=& (S_{D,W_s}S_{D,W_s}^T+\mu I)^{-1}+ \mu (S_{D,W_s}S_{D,W_s}^T+\mu I)^{-2} \nonumber\\
		&+&
		\mu^2  (S_{D,W_s}S_{D,W_s}^T+\mu I)^{-2}(S_{D,W_s}S_{D,W_s}^T)^{-1}.
\end{eqnarray}
Since for any  $h:[0,\infty)\rightarrow \mathbb R$ and positive operator $A$, there holds \cite{engl1996regularization}
\begin{equation}\label{exchange-1}
     h(A^TA) A^T=A^Th( AA^T) 
\end{equation}
for well defined $h(A^TA)$ and $h( AA^T) $.
it follows from \eqref{operator-KI-1},
 \eqref{Projection} and \eqref{def.noise-free}   that
\begin{eqnarray*}
  &&f_D^\diamond-f_D
  = 
  S_{D}^T(S_D S_{D}^T)^{-1}(S_Df^*-{\bf y}_D)\\
  &=&
  S_{D,W_s}^T(S_{D,W_s}S_{D,W_s}^T)^{-1}(S_{D,W_s}f^*-{\bf y}_{D,W_s})\\
  &=&
  S_{D,W_s}^T(S_{D,W_s}S_{D,W_s}^T+\mu I)^{-1} (S_{D,W_s}f^*-{\bf y}_{D,W_s})\\
  &+&
  \mu S_{D,W_s}^T  (S_{D,W_s}S_{D,W_s}^T+\mu I)^{-2}
  (S_{D,W_s}f^*-{\bf y}_{D,W_s})\\
  &+&
  \mu^2  S_{D,W_s}^T(S_{D,W_s}S_{D,W_s}^T)^{-1}(S_{D,W_s}S_{D,W_s}^T+\mu I)^{-2}  (S_{D,W_s}f^*-{\bf y}_{D,W_s})\\
  &=&
  (L_{\phi,D,W_s}+\mu I)^{-1}(L_{\phi,D,W_s}f^*-S_{D,W_s}^T{\bf y}_D)\\
  &+&
  \mu (L_{\phi,D,W_s}+\mu I)^{-2}(L_{\phi,D,W_s}f^*-S_{D,W_s}^T{\bf y}_D)\\
  &+&
  \mu^2(L_{\phi,D,W_s}+\mu I)^{-3/2}S_{D,W_s}^T(S_{D,W_s}S_{D,W_s}^T)^{-1}\\
  &&(S_{D,W_s}S_{D,W_s}^T+\mu I)^{-1/2}(S_{D,W_s}f^*-{\bf y}_{D,W_s}).
\end{eqnarray*}
Since for any  ${\bf c}=(c_1,\dots,c_{|D|})^T$, there holds
$$
    S_{D,W}^T(S_{D,W}S_{D,W}^T)^{-1}{\bf c}
    =S_D^T(S_DS_D^T)^{-1}\sqrt{n}W_s^{1/2}{\bf c}
    =S_D^T(S_DS_D^T)^{-2}S_DS_{D,W}^T{\bf c},
 $$
we obtain  from the Cordes inequality \eqref{Codes-inequality}, Lemma \ref{Lemma:integral-operator-relation} and    the inequality \cite{lin2020kernel} 
\begin{equation}\label{proj-operator-norm-1}
	\|S_D^T(S_DS_D^T)^{-2}S_D\|\leq  |D|(\sigma_{|D|}(\Phi_D))^{-1}
\end{equation}
 that for any $\mu\geq 16\tilde{c}'^2s^{-2\gamma}$, there holds
\begin{eqnarray*}
 &&\|f_D^\diamond-f_D\|_\psi
 \leq \|L_\phi^{\frac{1-\beta}2}(L_{\phi,D,W_s}+\mu I)^{-1}(L_{\phi,D,W_s}f^*-S_{D,W_s}^T{\bf y}_D)\|_\phi\\
 &+&
 \mu\|L_\phi^{\frac{1-\beta}2}  (L_{\phi,D,W_s}+\mu I)^{-2}(L_{\phi,D,W_s}f^*-S_{D,W_s}^T{\bf y}_D)\|_\phi\\
 &+&
 \mu^2\|L_\phi^{\frac{1-\beta}2} (L_{\phi,D,W_s}+\mu I)^{-3/2}S_D^T(S_DS_D^T)^{-2}S_D\|\\
 &\times&
 \| (L_{\phi,D,W_s}+\mu I)^{-1/2}(L_{\phi,D,W_s}f^*-S_{D,W_s}^T{\bf y}_{D,W_s})\|_\phi\\
 &\leq&
 4\mu^{-\frac\beta2} \left\|(  L_\phi+\mu I)^{-1/2}(L_{\phi,D,W_s}f^*-S_{D,W_s}^T{\bf y}_D)\right\|_\phi\\
 &+&
 2\mu^{1-\beta/2}|D|(\sigma_{|D|}(\Phi_D))^{-1}   \left\|(  L_\phi+\mu I)^{-1/2} (L_{\phi,D,W_s}f^*-S_{D}^T{\bf y}_{D,W_s})\right\|_\phi.
\end{eqnarray*}
Therefore,  Lemma \ref{lemma:value-difference-random} implies that for any $ \mu\geq 16\tilde{c}'^2s^{-2\gamma}$, with confidence $1-\delta$, there holds
$$
  \|f_D^\diamond-f_D\|_\psi\leq
   \mu^{-\frac\beta2}  (1+\mu |D| (\sigma_{|D|}(\Phi_D))^{-1}) 
  \bar{c}_3M\mu^{-\frac{d}{4\gamma}} |D|^{-1/2}\log\frac3\delta. 
$$
 This completes the proof of Lemma \ref{Lemma:Sample-error} by noting $s\sim |D|^{1/d}$.
\end{proof}

 Based on  Lemma \ref{Lemma:Approx-error} and Lemma \ref{Lemma:Sample-error}, we can derive the following lemma.

\begin{lemma}\label{Lemma:KI-random-noise}
Suppose that $\phi$, $\psi$ and $\varphi$ are given as  in  Lemma \ref{Lemma:Narcowich-interpolation} and that $f_D$ the 
minimal $\phi$-norm kernel interpolant defined in \eqref{minimal-norm-interpolation} based on noisy data (of the type \eqref{Model1:fixed}). Then there are two constants $C_1$ and $C_2$ depending only on $\alpha,\beta,\gamma,\tau$ and $d$ such that for any $f\in \mathcal N_\varphi$,
$0<\delta<1$,  and $\mu\geq C_1|D|^{-2\gamma/d}$, there holds
\begin{align}
\label{app-KI-random}
 \|f^*-f_D\|_{\psi}  \leq &
  C_2 \Bigl\{|D|^{\frac{(\beta-\alpha)\gamma}d}\|f^*\|_\varphi\\
& +  M \mu^{-\frac{2\gamma\beta\ + d}{4\gamma}}  \left[1+\mu |D| (\sigma_{|D|}(\Phi_D))^{-1}\right] 
   |D|^{-1/2}\log\frac3\delta\Bigr\} \nonumber
 \end{align}
with confidence at least $1-\delta$.
\end{lemma}

\begin{proof} 

Plugging \eqref{sample.error-1-fixed} and \eqref{approx.error-1} into \eqref{first-error-dec}, with confidence $1-\delta$, there holds
$$
 \|f^*-f_D\|_{\psi}  \leq 
  C_4 (\lambda^{\frac{\alpha-\beta}2}\|f^*\|_\varphi+  \mu^{-\frac\beta2}  (1+\mu |D| (\sigma_{|D|}(\Phi_D))^{-1}) 
  M\mu^{-\frac{d}{4\gamma}}|D|^{-1/2}\log\frac3\delta   
$$ 
for any $\mu,\lambda\geq C_3|D|^{-\frac{2\gamma}{d}}$.  Then
\eqref{app-KI-random} follows from $\lambda=C_3|D|^{-\frac{2\gamma}{d}}$.
 This completes the proof of Lemma \ref{Lemma:KI-random-noise}. 
\end{proof}

The following lemma provided in \cite[Example 2.10]{narcowich1998stability} presents a lower bound of  $\sigma_{|D|}(\Phi_D)$.

\begin{lemma}\label{Lemma:condition-number}
If $\hat{\phi}_k\sim k^{-2\gamma}$ with $2\gamma>d$ an integer, then
\begin{equation}\label{condition-number}
    \sigma_{|D|}( \Phi_D)\geq C'q_\Lambda^{4\gamma-d+1},
\end{equation}
where   $C'$  is a constant depending only on $d,\gamma$.
\end{lemma}

To prove Theorem \ref{Theorem:DKI-random}, we also need the following lemma that can be derived by using the same approach as \cite{lin2021distributed}.

\begin{lemma}\label{Lemma:distributed.1}
 For $\overline{f}_D$ defined by \eqref{DKI}, there  holds 
$$
         \mathbf E[\|\overline{f}_D-f^*\|_{\psi}^2] 
         \quad\leq
         \sum_{j=1}^m\frac{|D_j|^2}{|D|^2}\mathbf E[\|f_{D_j}
         -f^*\|_{\psi}^2]
        +\sum_{j=1}^m\frac{|D_j|}{|D|}
        \left\| f^\diamond_{D_j}-f^*\right\|_{\psi}^2,\label{eq:distributed.bound}
$$
where $f^{\diamond}_{D_j}$ is given by \eqref{Projection}.
\end{lemma}

We then use the above lemmas to prove Theorem \ref{Theorem:DKI-random}. 
\begin{proof}[Proof of Theorem \ref{Theorem:DKI-random}]
Due to Lemma \ref{Lemma:Approx-error}, we have 
\begin{equation}\label{approx.error-1-distributed}
 \|f_{D_j}^\diamond-f^*\|_\psi^2
    \leq \bar{C}_1|D_j|^{\frac{2\beta-2\alpha}d}\|f^*\|^2_\varphi
\end{equation}
where $\bar{C}_1$ is a constant depending only on $\tilde{c}$, $c$, $\alpha$,    $\beta$ and $\gamma$.
Furthermore, it follows from Lemma \ref{Lemma:KI-random-noise} that for any $\mu_j\geq C_3|D_j|^{-\frac{2\gamma}{d}}$, with confidence $1-\delta$, there holds
\begin{align*}
   \|f^*-f_{D_j}\|_{\psi} \leq &
  C_4 (|D_j|^{\frac{\beta-\alpha}d}\|f^*\|_\varphi\\
  &+  \mu_j^{-\frac\beta2}  (1+\mu_j |D_j| (\sigma_{|D_j|}(\Phi_{D_j}))^{-1}) 
   \mu_j^{-\frac{d}{4\gamma}}|D_j|^{-\frac12}M\log\frac3\delta.   
\end{align*}

Therefore, a direct calculation yields
\begin{align*}
  \mathbf E[\|f^*-f_{D_j}\|_{\psi}^2]
  \leq &
  \bar{C}_2 (|D_j|^{\frac{2\beta-2\alpha}d}\|f^*\|^2_\varphi\\
  &+  \mu_j^{-\beta}  (1+\mu_j |D_j| (\sigma_{|D_j|}(\Phi_{D_j}))^{-1})^2 
   \mu_j^{-\frac{d}{2\gamma}}|D_j|^{-1}M^2, 
\end{align*}
where $\bar{C}_2$ is a  constant depending only on $\tilde{c}$, $c$, $\alpha$,    $\beta$ and $\gamma$. Recalling Lemma \ref{Lemma:condition-number} and setting $\mu_j= C_3|D_j|^{-\frac{2\gamma}{d}}$, we have from the above inequality that
\begin{equation}\label{app-KI-random-distributed}
  \mathbf E[\|f^*-f_{D_j}\|_{\psi}^2]
  \leq
  \bar{C}_3 (|D_j|^{\frac{2\beta-2\alpha}d}\|f^*\|^2_\varphi+  |D_j|^{\frac{4\gamma+2\gamma\beta+2}d}   M^2),
\end{equation}
where $\bar{C}_3$ is a  constant depending only on $\tilde{c}$, $c$, $\alpha$,    $\beta$ and $\gamma$.
Plugging \eqref{app-KI-random-distributed} and \eqref{approx.error-1-distributed} into Lemma \ref{Lemma:distributed.1}, we then get
\begin{eqnarray*}
         \mathbf E[\|\overline{f}_D-f^*\|_{\psi}^2] 
         &\leq&
         \bar{C}_3 \sum_{j=1}^m\frac{|D_j|^2}{|D|^2}    |D_j|^{\frac{4\gamma+2\gamma\beta+2}d}   M^2 
         + (\bar{C}_2+1)\sum_{j=1}^m\frac{|D_j|}{|D|}
         |D_j|^{\frac{2\beta-2\alpha}d}\|f^*\|^2_\varphi.
\end{eqnarray*}
 This  completes the proof of Theorem \ref{Theorem:DKI-random}
\end{proof}

\section{Numerical Verifications}\label{Sec.Numerical}

Four simulations are carried out in this section to verify the excellent performance of DKI. The first one shows that DKI succeeds in circumventing the uncertainty of kernel interpolation. The second one exhibits the role of $m$ in DKI. The third one focuses on the role of the division strategy in DKI. The last one compares DKI with several popular spherical data fitting schemes including the  distributed filtered hyperinterpolation (DFH) \cite{lin2021distributed}, sketching with $s^*$-designs \cite{lin2021subsampling}, and distributed kernel ridge regression (DKRR) \cite{Feng2021radial}.


The training  samples are generated according to \ref{Model1:fixed} with $\varepsilon_i$  i.i.d. being drawing i.i.d. according to the  Gaussian distribution $\mathcal{N}(0, \delta^2)$ and  the target function
\begin{equation}\label{f1}
f({x}) = \sum\limits_{i=1}^{\kappa} \tilde{\psi}\left(\frac{\|{x}-{z}_i\|_2}{c}\right),
\end{equation}
where $\tilde{\psi}(u)$ is the well-known Wendland function \cite{chernih2014wendland}
\begin{equation}\label{psi_wave}
\tilde{\psi}(u) = (1-u)_{+}^8(32u^3+25u^2+8u+1),
\end{equation} 
$u_+ = \max\{u, 0\}$, and ${z}_i \in \mathbb S^{d}$ ($i=1,\cdots,\kappa$) are the center points of the regions of an equal area partitioned by Leopardi's recursive zonal sphere partitioning procedure \cite{leopardi2006partition}\footnote{http://eqsp.sourceforge.net}. 
The inputs $\{{x}_i\}_{i=1}^N$ of training samples are generated by Womersley's symmetric spherical $t$-designs with $N=(t)^2/2+t/2+O(1)$ points \cite{womersley2018efficient} \footnote{https://web.maths.unsw.edu.au/\%7Ersw/Sphere/EffSphDes/}. The testing samples are generated according to the clean model $y_j'=f(x_j')$ whose inputs 
 $\{x_j'\}_{j=1}^{N'}$  are $N'=10000$ equally distributed spiral points on $\mathbb S^2$ generated by $\left[\sin\alpha_j\cos\beta_j,\sin\alpha_j\sin\beta_j,\cos\alpha_j)\right]^\top$ with $\alpha_j=\arccos(1-(2j-1)/N')$ and $\beta_j=\mod(1.8\sqrt{N'}\theta_j, 2\pi)$.   We take $\phi_1(x_1, x_2)=\tilde{\psi}(\|x_1 - x_2\|_2)$ as the positive definite kernel. All parameters in the simulations are selected by grid search. Each simulation is conducted 30 times, and the average results are recorded.

\begin{figure*}[t]
    \centering
     \subfigure[Uncertainty of KI]{\includegraphics[width=6cm,height=4.5cm]{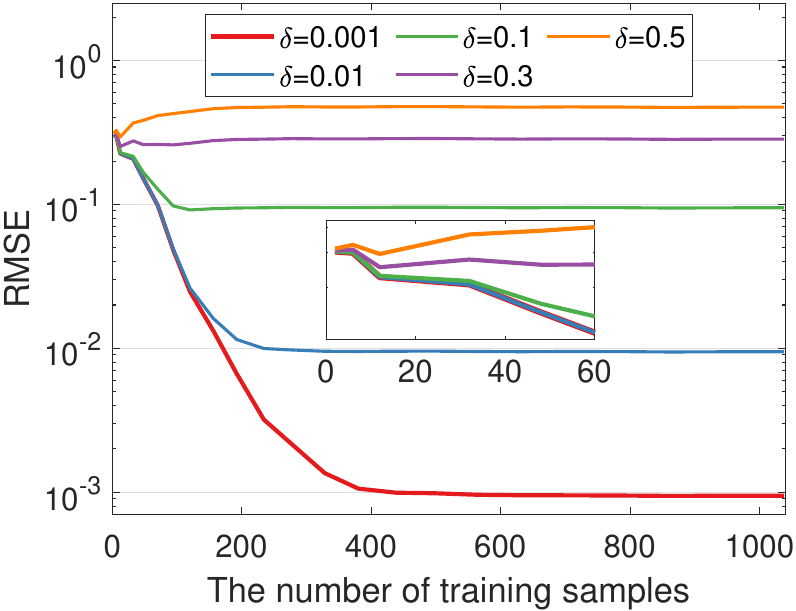}}
    \subfigure[Overcoming uncertainty by DKI ]{\includegraphics[width=6cm,height=4.5cm]{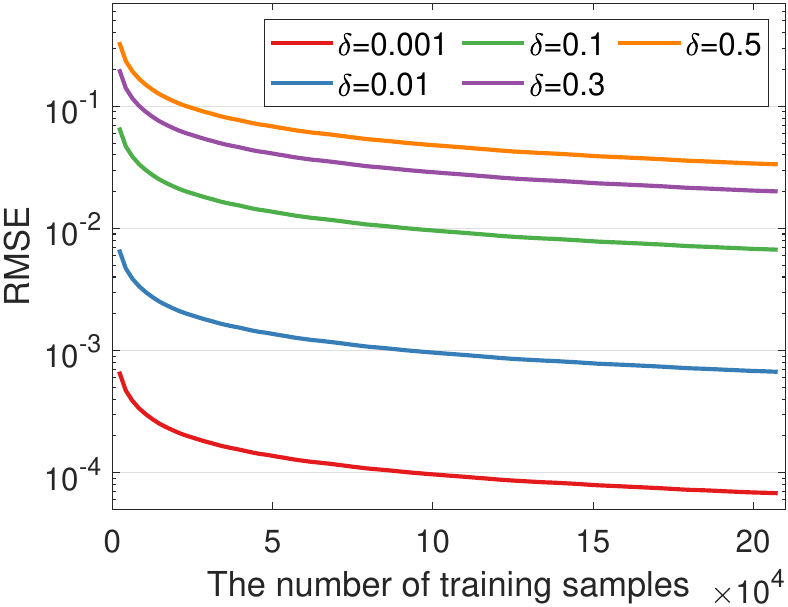}}
    \caption{Relations between the RMSE and the number of training samples for KI and DKI under different levels of Gaussian noise}\label{sim1_RMSE_increase_sam}
\end{figure*}

\textbf{Simulation 1:} In this simulation, we show that DKI can circumvent the uncertainty of KI via showing the relations between the RMSE and data size for different noise levels. We firstly exhibit the uncertainty of KI. $t$-designs  with varying $t\in\{1,3,\cdots,45\}$ are generated  and  $\delta\in\{0.001,0.01,0.1,0.3,0.5\}$ are taken as the training samples and noise level, respectively . The curves of the testing RMSE of kernel interpolation changing with the number of training samples are plotted in Figure \ref{sim1_RMSE_increase_sam} (a).
It can be seen that as the number of training samples grows, the RMSE is decreasing for small noise (i.e., $\delta\leq 0.1$),  and non-decreasing  for high noises (i.e., $\delta>0.1$), showing the uncertainty of KI.


To show the outperformance of DKI, we involve more samples.  The inputs of the  training data is generated by the following two steps: First, Womersley's symmetric spherical $45$-designs are used to generate $1038$ points $\{\hat{x}_{i}\}_{i=1}^{1038}$. Secondly, the points $\{\hat{x}_{i}\}_{i=1}^{1038}$ are rotated by the rotation matrix
$$
A_k := \left(
  \begin{array}{ccc}
    \cos(k\pi/10) & -\sin(k\pi/10) & 0 \\
    \sin(k\pi/10) & \cos(k\pi/10) & 0 \\
    0 & 0 & 1 \\
  \end{array}
\right)
$$
to obtain the new points $\{{x}_{i,k}\}_{i=1}^{1038}$ for $k=1,\cdots,200$, i.e., ${x}_{i,k}=A_k\hat{x}_{i}$, and the points $\{x_{i,k}\}_{i=1}^{1038}$ are used as the inputs of training samples on the $k$th local machine. That is, there are totally  $1038k$ samples for $k=2,4,\dots,200$. 
The curves of the testing RMSE of DKI changing with $k$ are shown in Figure \ref{sim1_RMSE_increase_sam} (b). From the results, we can see that the RMSE monotonously decreases with respect to the size of samples for all levels of Gaussian noise, which is consistent with our assertion that more samples would yield a higher-quality approximation. This also demonstrates that the proposed divide-and-conquer approach can effectively handle the uncertainty of kernel interpolation.


\begin{figure*}[t]
    \centering
     {\includegraphics[width=7cm,height=5cm]{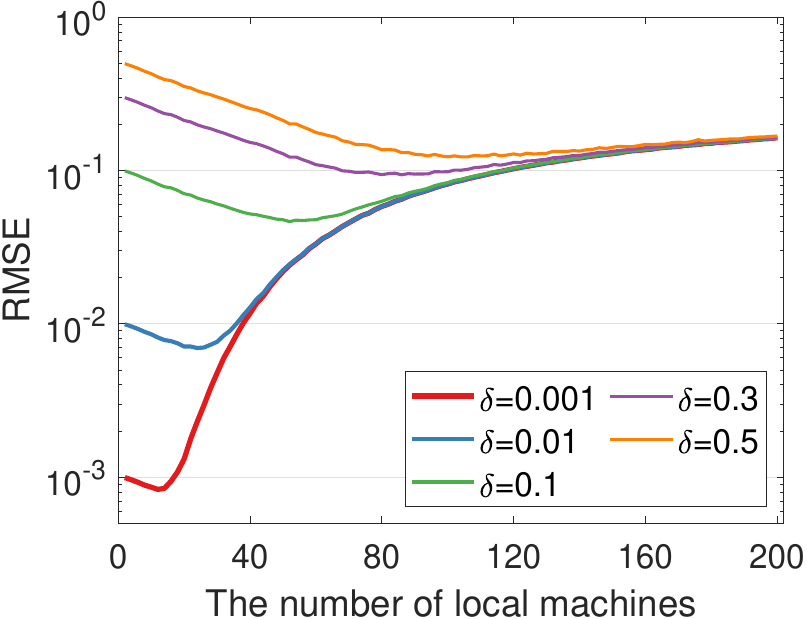}}
    \caption{Relation between   RMSE and the number of local estimators for DKI under different levels of Gaussian noise}\label{Sim2_RMSE}
\end{figure*}

\textbf{Simulation 2:} In this simulation,  we show the role of the parameter $m$ in DKI.
We generate $10014$ training samples (with 141-designs as inputs).  
The number  of divisions, $m$, ranges from $\{5,10,\cdots,200\}$. Figure \ref{Sim2_RMSE} shows the relation between RMSE of DKI and the number of local machines under different levels of Gaussian noise, provided that the total number of training samples is given.  From Figure \ref{Sim2_RMSE}, we can conclude the following assertions: 1) For training samples with higher levels of noise, the testing RMSE generally decreases at first and then increases slowly as the number of local machines increases. Moderate
values of $m$ are more conductive to good approximation property for DKI. The reason is that too small $m$  does not successfully address the uncertainty issue in kernel interpolation; too large $m$ increases the fitting error, resulting in slightly worse generalization performance. 2) The optimal number $m$ with the lowest RMSE grows with increasing Gaussian noise.
This verifies the equation \eqref{dist-error} of Theorem \ref{Theorem:DKI-random}, in which the approximation error is primarily concerned with the sample error for large noise (i.e., large $M$) and can be reduced using a large $m$. 



\begin{figure*}[t]
    \centering
    {\includegraphics[width=7cm,height=3cm]{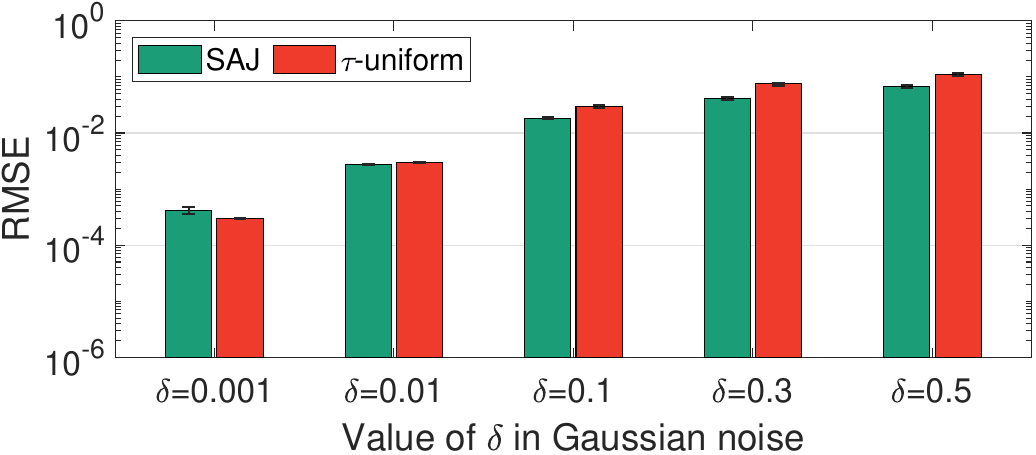}}
    \caption{Comparison of the testing RMSE for $\tau$-uniform and SAJ data divisions with increasing levels of Gaussian noise on the 3-dimensional data. }\label{DivideComparison}
\end{figure*}

\begin{figure*}[t]
    \centering
    \subfigure[RMSE]{\includegraphics[width=6.5cm,height=3.5cm]{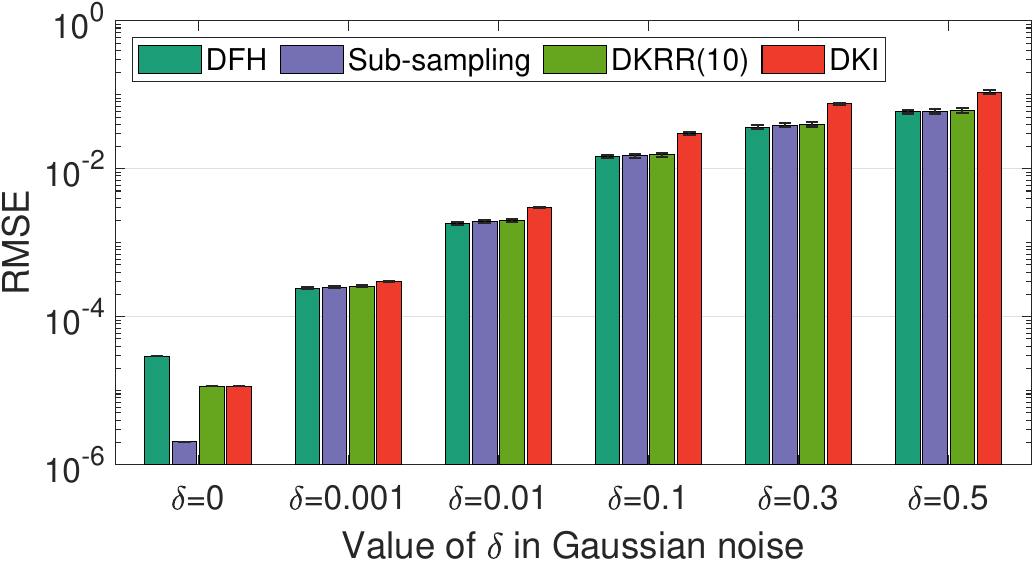}}\hspace{-0.04in}
    \subfigure[Training time]{\includegraphics[width=6.5cm,height=3.5cm]{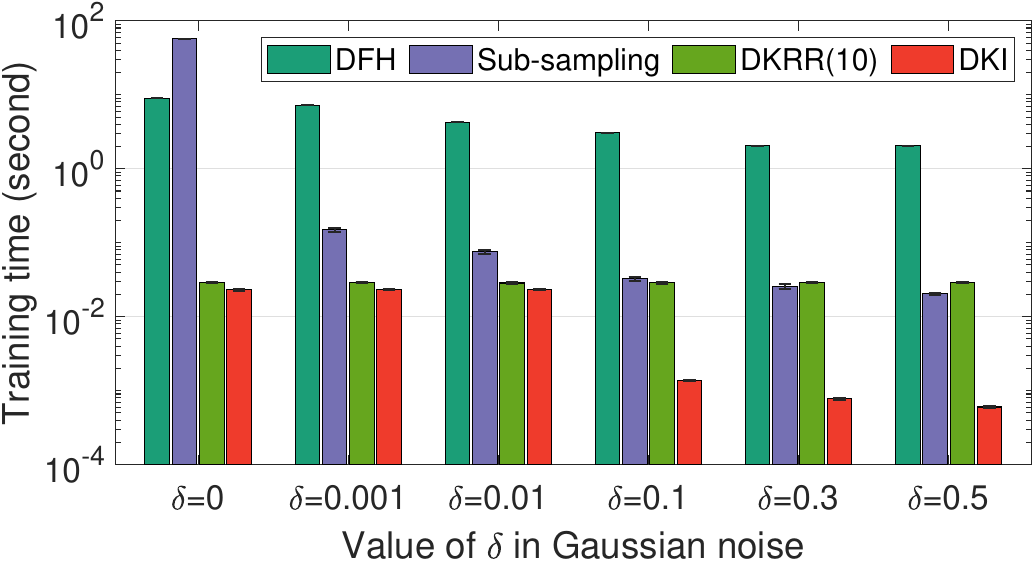}}
    \caption{The comparison of RMSE and training time among the four methods with increasing levels of Gaussian noise on the 3-dimensional data.}\label{Sim5_dim3}
\end{figure*}

\begin{figure*}[htbp]
    \centering
    \subfigcapskip=-2pt
    \subfigure[Groudtruth]{\includegraphics[width=3.1cm, height=3.7cm]{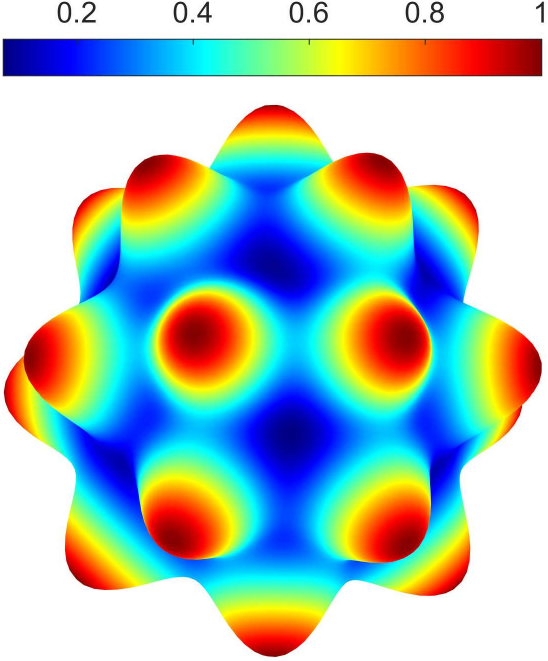}}
    \subfigure[$\delta=0.001$]{\includegraphics[width=3.1cm, height=3.7cm]{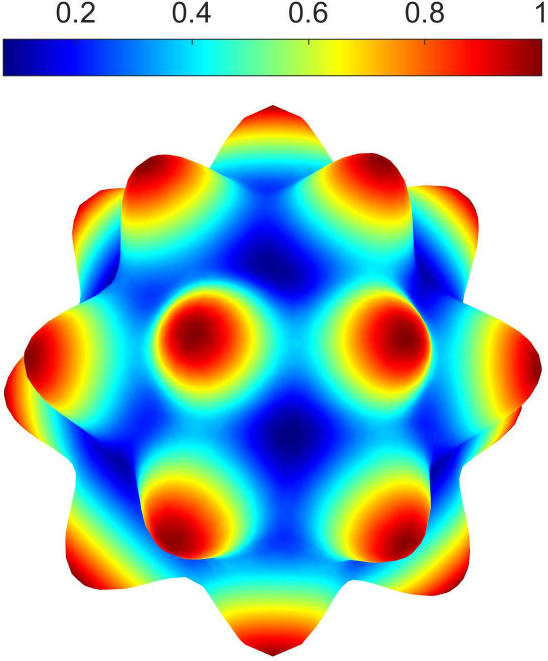}}
    \subfigure[$\delta=0.1$]{\includegraphics[width=3.1cm, height=3.7cm]{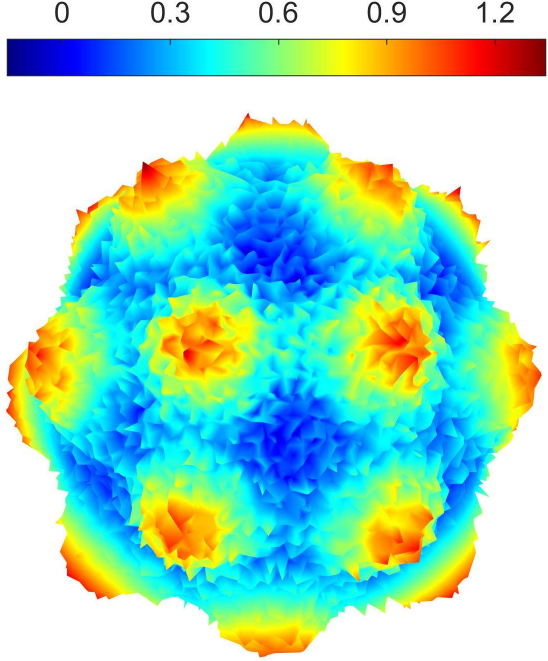}}
    \subfigure[$\delta=0.5$]{\raisebox{0.055\height}{\includegraphics[width=3.1cm, height=3.5cm]{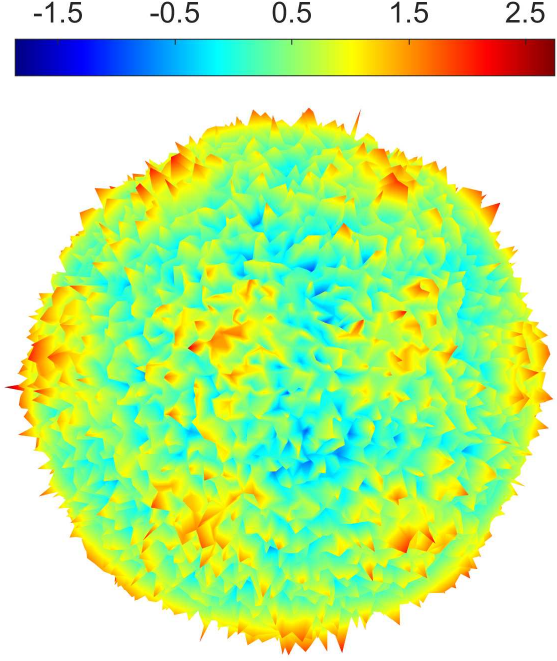}}}
	\caption{Visualization of the 3-dimensional training data with different levels of Gaussian noise.}\label{GroudtruthAddNoise}
\end{figure*}

\begin{figure}[htbp]
	\centering
	\subfigure{
        \rotatebox{90}{\small{~~~~~~~~~$\delta=0.001$}}
    	\includegraphics[width=3cm, height=3.5cm]{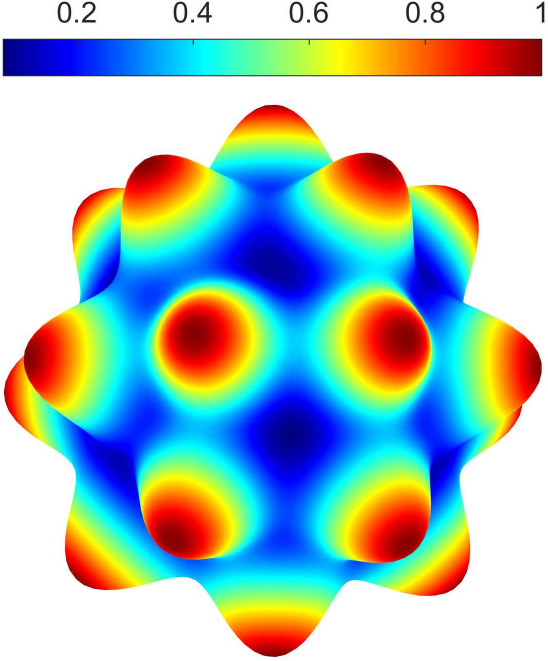}
	}\hspace{-0.1in}
		\subfigure{
		\includegraphics[width=3cm, height=3.5cm]{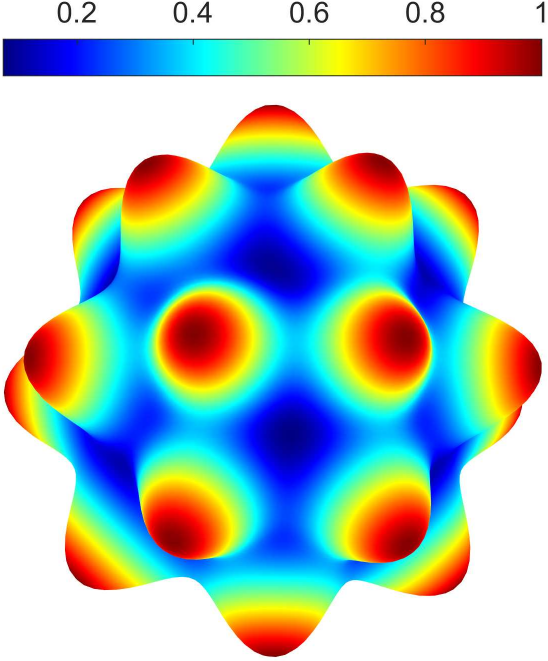}
	}\hspace{-0.1in}
		\subfigure{
		\includegraphics[width=3cm, height=3.5cm]{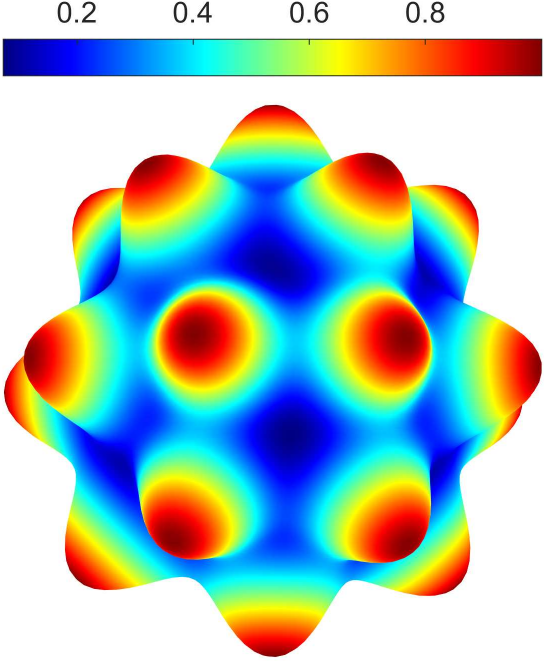}
	}\hspace{-0.1in}
		\subfigure{
	    \includegraphics[width=3cm, height=3.5cm]{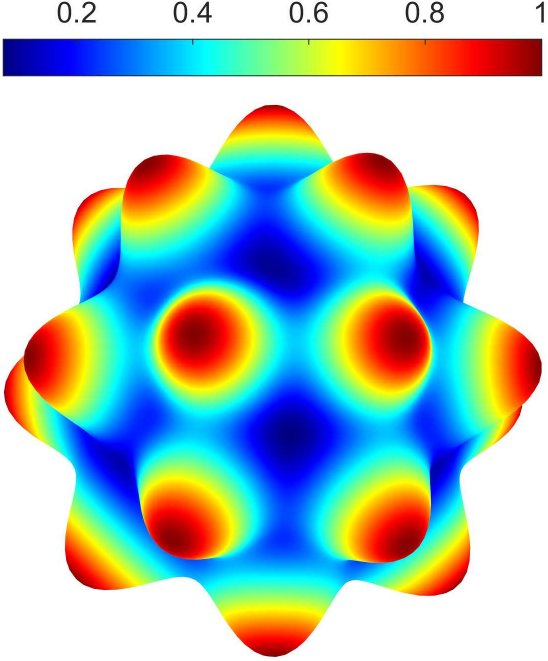}
	}
	\vspace{-3mm}
		\subfigure{
        \rotatebox{90}{\small{~~~~~~~~~$\delta=0.1$}}
		\includegraphics[width=3cm, height=3.5cm]{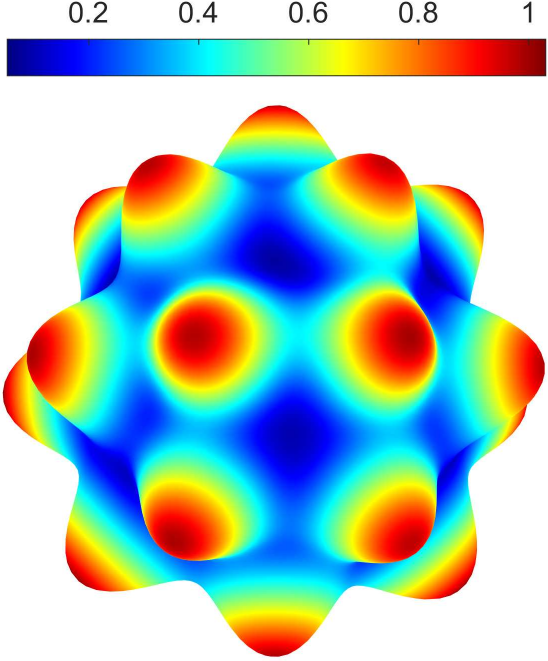}
	}\hspace{-0.1in}
		\subfigure{
		\includegraphics[width=3cm, height=3.5cm]{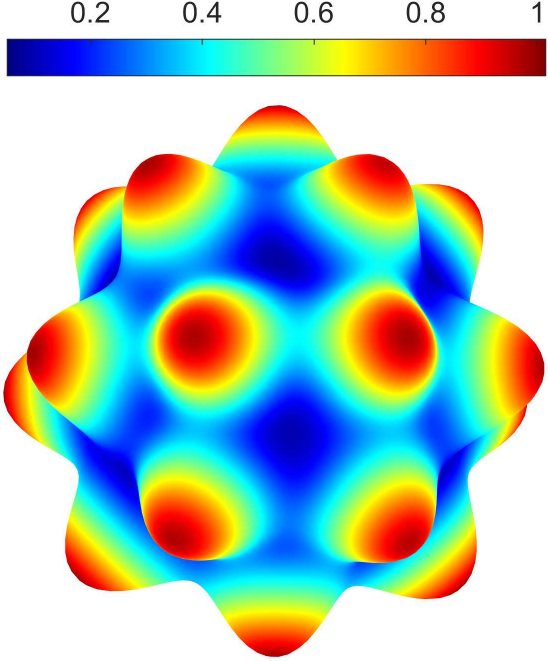}
	}\hspace{-0.1in}
		\subfigure{
		\includegraphics[width=3cm, height=3.5cm]{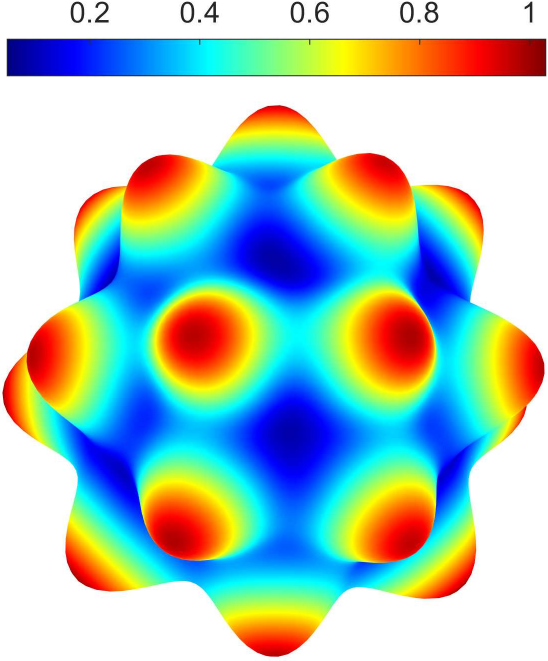}
	}\hspace{-0.1in}
		\subfigure{
		\includegraphics[width=3cm, height=3.5cm]{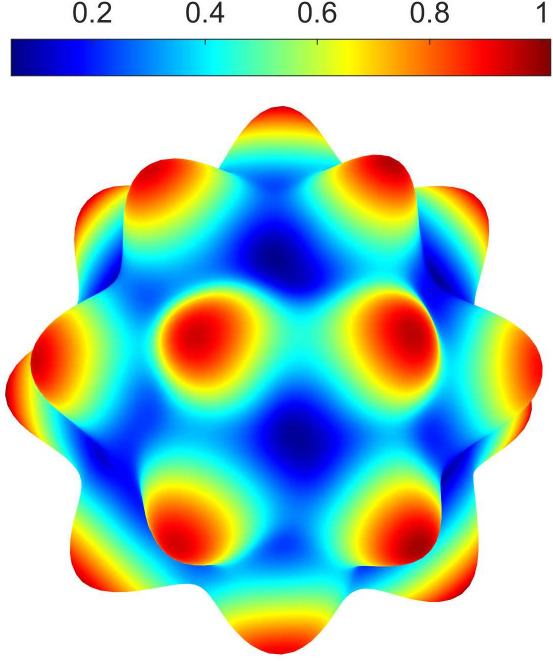}
	}
	\vspace{-3mm}
	\setcounter{subfigure}{0}
    \subfigure[DFH]{
		\rotatebox{90}{\small{~~~~~~~~~$\delta=0.5$}}
		\includegraphics[width=3cm, height=3.5cm]{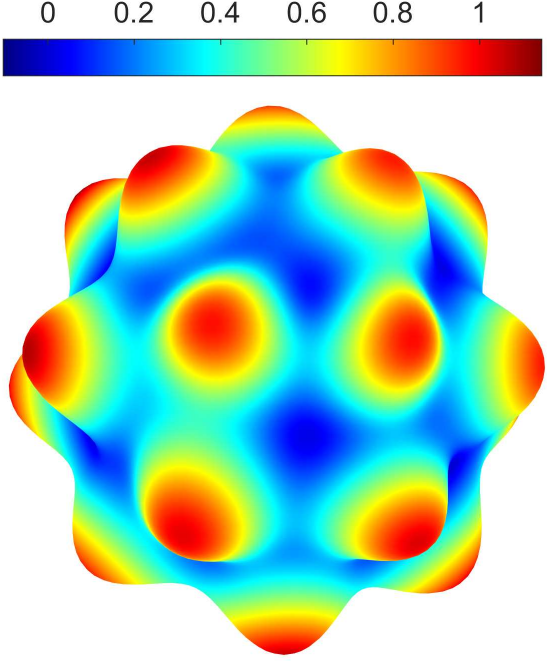}
	}\hspace{-0.11in}
	\subfigure[Subsampling]{
		\includegraphics[width=3cm, height=3.5cm]{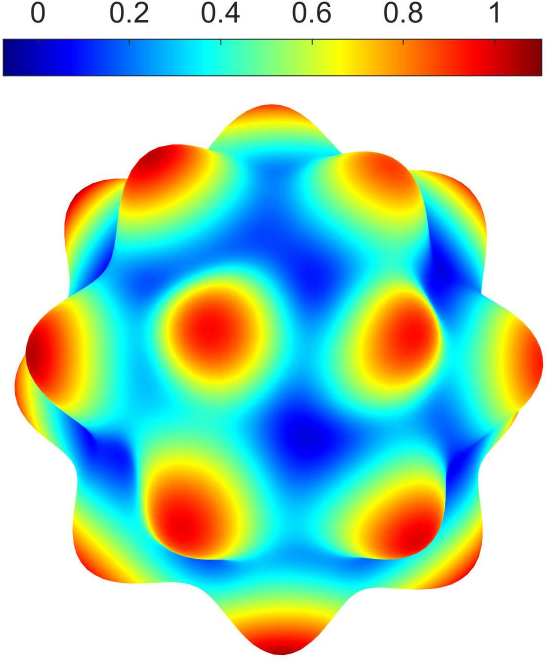}
	}\hspace{-0.11in}
		\subfigure[DKRR(10)]{
		\includegraphics[width=3cm, height=3.5cm]{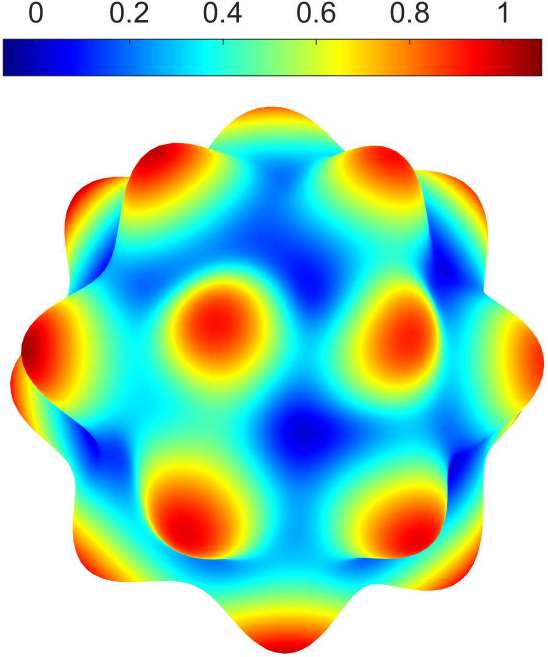}
	}\hspace{-0.11in}
		\subfigure[DKI]{
		\includegraphics[width=3cm, height=3.5cm]{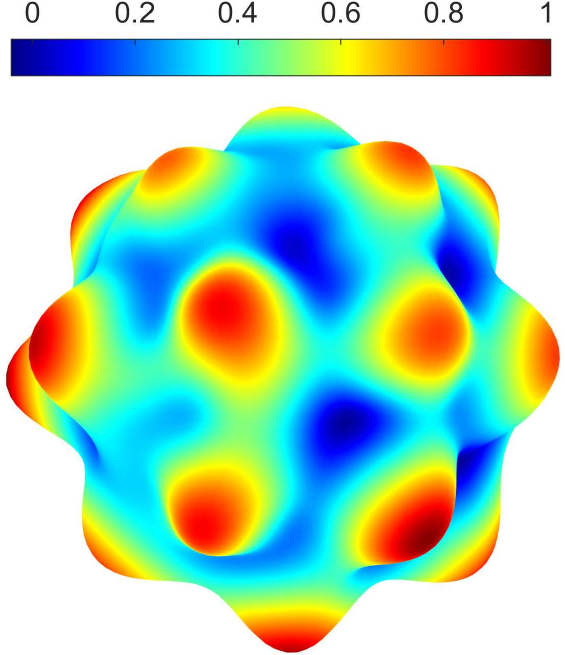}
	}
	\caption{Visualization of recovered results on the 3-dimensional data with different noise standard deviations.}
	\label{VisualRecovery}
\end{figure}

\textbf{Simulation 3:} This simulation study of the role of the division strategy that is used in Step 1 of DKI. The data generation is the same as in Simulation 2, but with 10 rotation matrices, so we generate 10 sets of samples $S_k=\{(x_{i,k}, y_{i,k})\}_{i=1}^{1038}$ for $k=1,\cdots,10$ with a total of $N=10380$ training samples. The details of data division based on 
$\tau$-uniform sets are as follows: Let $m$ ($m\geq 10$) be the number of local machines, and $r:=\mod(m,10)$. If $r=0$, that is, if $m$ can be divided by $10$, then the samples of each set $S_k$ are distributed randomly and equally to $m/10$ local machines for $k=1,\cdots,10$. If $r>0$, we randomly choose $r$ sets from $\{S_k\}_{k=1}^{10}$; samples from each chosen set are distributed equally to $\lceil m/10 \rceil$ local machines; samples of each set from the remaining $10-r$ sets are distributed equally to $\lfloor m/10 \rfloor$ local machines. The number of local machines for data division based on $\tau$-uniform is chosen from the set $\{10,12,\dots,100\}$. The radius value $c_0$ of SAJ is chosen from the set $\{0.05,0.1,\cdots,1\}$. The testing RMSE of the two data division methods with optimal parameters is shown in Figure \ref{DivideComparison}, from which we know that the two data division methods are comparable. Data division based on $\tau$-uniform sets is slightly worse than SAJ for large noise. This may be due to the fact that the optimal number of local machines is typically large for large noise, resulting in some differences in the number of samples distributed to the local machines using data division based on $\tau$-uniform sets (e.g., some local machines have $1038/\lceil m/10 \rceil$ training samples and other machines have $1038/\lfloor m/10 \rfloor$ training samples), whereas the sample number equalization steps in SAJ make the sample number in each local machine as close as possible. Hence, the proposed SAJ approach is a feasible scheme for dividing the spherically scattered data into disjoint $\tau$-quasi uniform sets.


\textbf{Simulation 4:}  In this simulation, we   show the efficiency and effectiveness of the proposed DKI when compared with DFH, sketching with $s^*$-design, and DKRR. To satisfy the training requirements of DFH, the data generation is the same as in Simulation 3. For DFH, the training samples are equally distributed to $10$ local machines, i.e., the samples $S_k:=\{(x_{i,k}, y_{i,k})\}_{i=1}^{1038}$ obtained by the $k$th rotation matrix are located on the $k$th local machine, and the parameter $L$ related to the polynomial degree is selected from the set $\{2,4,\dots,40\}$. For DKRR, as the number of local machines increases, the training complexity is reduced while the generalization performance is degraded. Hence, we make a trade-off between the training complexity and generalization performance and set the number of local machines to $10$, the same as DFH. The regularization parameter $\lambda$ of DKRR is chosen from the set 
$\{\frac{1}{2^q}:\frac{1}{2^q}>10^{-10},q=0,1,2,\dots\}$. For sub-sampling with $s^*$-designs, we choose $s^*$ from the set $\{1,3,\dots,121\}$. For DKI, we employ the data division based on $\tau$-uniform sets for simplicity and choose the number of local machines from the set $\{10,12,\dots,100\}$. Figure \ref{Sim5_dim3} records the testing RMSE and corresponding training time of the involved methods with the best parameters. The number in parentheses after `DKRR' is the number of local machines. In order to have an intuitive comparison, the 3-dimensional training data and the recovered results of the compared methods under different levels of Gaussian noise are visualized in Figures \ref{GroudtruthAddNoise} and \ref{VisualRecovery}, respectively. From the results, it can be seen that DKI is
highly efficient in the training stage while generalizing comparably to state-of-the-art methods, which implies that DKI is an effective approach to handling the uncertainty problem of kernel interpolation on spheres.

\bibliographystyle{siamplain}
\bibliography{distributed-interpolation}

\begin{thebibliography}{10}

\bibitem{bhatia2013matrix}
{\sc R.~Bhatia}, {\em Matrix Analysis}, vol.~169, Springer Science \& Business
  Media, 2013.

\bibitem{brauchart2007numerical}
{\sc J.~S. Brauchart and K.~Hesse}, {\em Numerical integration over spheres of
  arbitrary dimension}, Constructive Approximation, 25 (2007), pp.~41--71.

\bibitem{brown2005approximation}
{\sc G.~Brown and F.~Dai}, {\em Approximation of smooth functions on compact
  two-point homogeneous spaces}, Journal of Functional Analysis, 220 (2005),
  pp.~401--423.

\bibitem{chernih2014wendland}
{\sc A.~Chernih, I.~H. Sloan, and R.~S. Womersley}, {\em Wendland functions
  with increasing smoothness converge to a gaussian}, Advances in Computational
  Mathematics, 40 (2014), pp.~185--200.

\bibitem{dai2006multivariate}
{\sc F.~Dai}, {\em Multivariate polynomial inequalities with respect to
  doubling weights and $a_\infty$ weights}, Journal of Functional Analysis, 235
  (2006), pp.~137--170.

\bibitem{dai2006generalized}
{\sc F.~Dai}, {\em On generalized hyperinterpolation on the sphere},
  Proceedings of the American Mathematical Society, 134 (2006), pp.~2931--2941.

\bibitem{engl1996regularization}
{\sc H.~W. Engl, M.~Hanke, and A.~Neubauer}, {\em Regularization of Inverse
  Problems}, vol.~375, Springer Science \& Business Media, 1996.

\bibitem{Feng2021radial}
{\sc H.~Feng, S.-B. Lin, and D.-X. Zhou}, {\em Radial basis function
  approximation with distributively stored data on spheres}, arXiv:2112.02499,
  (2021).

\bibitem{hangelbroek2011kernel}
{\sc T.~Hangelbroek, F.~J. Narcowich, X.~Sun, and J.~D. Ward}, {\em Kernel
  approximation on manifolds ii: The $l_\infty$ norm of the $l_2$ projector},
  SIAM Journal on Mathematical Analysis, 43 (2011), pp.~662--684.

\bibitem{hangelbroek2010kernel}
{\sc T.~Hangelbroek, F.~J. Narcowich, and J.~D. Ward}, {\em Kernel
  approximation on manifolds i: bounding the lebesgue constant}, SIAM Journal
  on Mathematical Analysis, 42 (2010), pp.~1732--1760.

\bibitem{hesse2017radial}
{\sc K.~Hesse, I.~H. Sloan, and R.~S. Womersley}, {\em Radial basis function
  approximation of noisy scattered data on the sphere}, Numerische Mathematik,
  137 (2017), pp.~579--605.

\bibitem{hesse2021local}
{\sc K.~Hesse, I.~H. Sloan, and R.~S. Womersley}, {\em Local rbf-based
  penalized least-squares approximation on the sphere with noisy scattered
  data}, Journal of Computational and Applied Mathematics, 382 (2021),
  p.~113061.

\bibitem{hubbert2015spherical}
{\sc S.~Hubbert, Q.~T. L{\^e}~Gia, and T.~M. Morton}, {\em Spherical radial
  basis functions, theory and applications}, Springer, 2015.

\bibitem{king2012lower}
{\sc M.~A. King, R.~J. Bingham, P.~Moore, P.~L. Whitehouse, M.~J. Bentley, and
  G.~A. Milne}, {\em Lower satellite-gravimetry estimates of antarctic
  sea-level contribution}, Nature, 491 (2012), pp.~586--589.

\bibitem{le2006continuous}
{\sc Q.~T. Le~Gia, F.~J. Narcowich, J.~D. Ward, and H.~Wendland}, {\em
  Continuous and discrete least-squares approximation by radial basis functions
  on spheres}, Journal of Approximation Theory, 143 (2006), pp.~124--133.

\bibitem{leopardi2006partition}
{\sc P.~Leopardi}, {\em A partition of the unit sphere into regions of equal
  area and small diameter}, Electronic Transactions on Numerical Analysis, 25
  (2006), pp.~309--327.

\bibitem{levesley1999norm}
{\sc J.~Levesley, Z.~Luo, and X.~Sun}, {\em Norm estimates of interpolation
  matrices and their inverses associated with strictly positive definite
  functions}, Proceedings of the American Mathematical Society,  (1999),
  pp.~2127--2134.

\bibitem{lin2020kernel}
{\sc S.-B. Lin, X.~Chang, and X.~Sun}, {\em Kernel interpolation of high
  dimensional scattered data}, arXiv preprint arXiv:2009.01514,  (2020).

\bibitem{lin2017distributed}
{\sc S.-B. Lin, X.~Guo, and D.-X. Zhou}, {\em Distributed learning with
  regularized least squares}, The Journal of Machine Learning Research, 18
  (2017), pp.~3202--3232.

\bibitem{lin2021subsampling}
{\sc S.-B. Lin, D.~Wang, and D.-X. Zhou}, {\em Sketching with spherical designs
  on spheres}, SIAM Journal on Scientific Computation, in press (2023).

\bibitem{lin2021distributed}
{\sc S.-B. Lin, Y.~G. Wang, and D.-X. Zhou}, {\em Distributed filtered
  hyperinterpolation for noisy data on the sphere}, SIAM Journal on Numerical
  Analysis, 59 (2021), pp.~634--659.

\bibitem{mcewen2011novel}
{\sc J.~D. McEwen and Y.~Wiaux}, {\em A novel sampling theorem on the sphere},
  IEEE Transactions on Signal Processing, 59 (2011), pp.~5876--5887.

\bibitem{mhaskar2001spherical}
{\sc H.~Mhaskar, F.~Narcowich, and J.~Ward}, {\em Spherical
  marcinkiewicz-zygmund inequalities and positive quadrature}, Mathematics of
  Computation, 70 (2001), pp.~1113--1130.

\bibitem{mhaskar2006weighted}
{\sc H.~N. Mhaskar}, {\em Weighted quadrature formulas and approximation by
  zonal function networks on the sphere}, Journal of Complexity, 22 (2006),
  pp.~348--370.

\bibitem{muller2006spherical}
{\sc C.~M{\"u}ller}, {\em Spherical Harmonics}, vol.~17, Springer, 1966.

\bibitem{narcowich1998stability}
{\sc F.~J. Narcowich, N.~Sivakumar, and J.~D. Ward}, {\em Stability results for
  scattered-data interpolation on euclidean spheres}, Advances in Computational
  Mathematics, 8 (1998), pp.~137--163.

\bibitem{narcowich2007direct}
{\sc F.~J. Narcowich, X.~Sun, J.~D. Ward, and H.~Wendland}, {\em Direct and
  inverse sobolev error estimates for scattered data interpolation via
  spherical basis functions}, Foundations of Computational Mathematics, 7
  (2007), pp.~369--390.

\bibitem{narcowich2002scattered}
{\sc F.~J. Narcowich and J.~D. Ward}, {\em Scattered data interpolation on
  spheres: error estimates and locally supported basis functions}, SIAM Journal
  on Mathematical Analysis, 33 (2002), pp.~1393--1410.

\bibitem{rudi2015less}
{\sc A.~Rudi, R.~Camoriano, and L.~Rosasco}, {\em Less is more: Nystr{\"o}m
  computational regularization.}, in NIPS, 2015, pp.~1657--1665.

\bibitem{schaback1995error}
{\sc R.~Schaback}, {\em Error estimates and condition numbers for radial basis
  function interpolation}, Advances in Computational Mathematics, 3 (1995),
  pp.~251--264.

\bibitem{smale2004shannon}
{\sc S.~Smale and D.-X. Zhou}, {\em Shannon sampling and function
  reconstruction from point values}, Bulletin of the American Mathematical
  Society, 41 (2004), pp.~279--305.

\bibitem{smale2005shannon}
{\sc S.~Smale and D.-X. Zhou}, {\em Shannon sampling ii: Connections to
  learning theory}, Applied and Computational Harmonic Analysis, 19 (2005),
  pp.~285--302.

\bibitem{tsai2006all}
{\sc Y.-T. Tsai and Z.-C. Shih}, {\em All-frequency precomputed radiance
  transfer using spherical radial basis functions and clustered tensor
  approximation}, ACM Transactions on Graphics (TOG), 25 (2006), pp.~967--976.

\bibitem{wendland2004scattered}
{\sc H.~Wendland}, {\em Scattered data approximation}, vol.~17, Cambridge
  university press, 2004.

\bibitem{wieczorek1998potential}
{\sc M.~A. Wieczorek and R.~J. Phillips}, {\em Potential anomalies on a sphere:
  Applications to the thickness of the lunar crust}, Journal of Geophysical
  Research: Planets, 103 (1998), pp.~1715--1724.

\bibitem{womersley2018efficient}
{\sc R.~S. Womersley}, {\em Efficient spherical designs with good geometric
  properties}, in Contemporary computational mathematics-A celebration of the
  80th birthday of Ian Sloan, Springer, 2018, pp.~1243--1285.

\bibitem{zhang2015divide}
{\sc Y.~Zhang, J.~Duchi, and M.~Wainwright}, {\em Divide and conquer kernel
  ridge regression: A distributed algorithm with minimax optimal rates}, The
  Journal of Machine Learning Research, 16 (2015), pp.~3299--3340.

\end{thebibliography}

\end{document}


\maketitle

\section{Appendix A: Select-and-Judge Strategy for Data Division}
In this Appendix, we present the detailed implementation of the select-and-judge (SAJ) strategy. Our aim is to derive a series of subsets of similar Cardinality with   separation radius not smaller than a given tolerance $c_0$. There are two stages for SAJ.  

{\bf Stage 1: Division for the separation radius requirement}

$\bullet$ Step 1: Compute pairwise geodesic distances between samples in $\Lambda$.  

$\bullet$ Step 2: Select an arbitrary $x\in\Lambda$ and label it as $x_{1,j}$ for $j=1$.

$\bullet$ Step 3:  Define $\Lambda^1=\Lambda\backslash B(x_{1,j},c_0)$, where $B(x,c_0)$ the spherical cap with center $x_{1,j}$ and radius $c_0$. That is, to remove the points in $B(x_{1,j},c_0)$ from $\Lambda$.

$\bullet$ Step 4. Given a subset $\Lambda_j^\ell:=\{x_{1,j},\dots,x_{\ell,j}\}$, randomly select an  $x\in \Lambda^{\ell}$ and label it as $x_{\ell+1,j}$. Update $\Lambda_j^{\ell+1}=\Lambda_j^\ell\cup\{x_{\ell+1,j}\}$ and $\Lambda^{\ell+1}= {\Lambda^\ell}\backslash B(x_{\ell+1,j},c_0)$.

$\bullet$ Step 5. Repeat the above procedure until there is an
$U\in\mathbb N$ such that $\Lambda^{U}=\{\}.$ We then obtain
$\Lambda_j=\Lambda_{j}^{U-1}$ for $j=1$.

$\bullet$ Step 6: Repeat the above five steps with $\Lambda=\Lambda\backslash \Lambda_j$ and we then obtain $\Lambda_{j+1}$.

$\bullet$ Step 7: Stop the algorithm if all points are divided into data subsets and obtain $m$ disjoint subsets $\Lambda_1,\dots,\Lambda_m$ with $q_{\Lambda_j}>c_0$.

The numbers of samples in the obtained $m$ disjoint subsets $\Lambda_1,\dots,\Lambda_m$ are seriously unbalanced. In general, the   first subset contains more samples. This will degrade the performance of DKI. Hence, we design the following steps for sample number equalization, whose core idea is 
that some samples from the subsets with sample numbers greater than the average number are assigned to the subsets with sample numbers less than the average number.

{\bf Stage 2: Balance the size of data}

$\bullet$ Step 1: Calculate the average number of samples $\bar{N}$ in each local machine using $\bar{N} = \lfloor N/m \rfloor$.  

$\bullet$ Step 2: Divide the subsets $\Lambda_1,\dots,\Lambda_m$ into two parts: those with sample numbers greater than $\bar{N}$ and those with sample numbers not greater than $\bar{N}$. We assume, without sacrificing generality, that the sample numbers of the first $l$ subsets  $\Lambda_1,\dots,\Lambda_l$ are greater than $\bar{N}$ and that the sample numbers of the rest subsets $\Lambda_{l+1},\dots,\Lambda_m$ are not greater than $\bar{N}$.


$\bullet$ Step 3:  For $j=1,\cdots,l$, randomly divide $\Lambda_j$ into two disjoint subsets, denoted as $\Lambda_j^1$ and $\Lambda_j^2$, where $|\Lambda_j^1|=\bar{N}$ and $|\Lambda_j^2|=|\Lambda_j|-|\Lambda_j^1|$.

$\bullet$ Step 4. Select an arbitrary $x\in\Lambda_j^2$, if there exists $j_0>l$ such that $q_{\Lambda_{j_0}\cup x}>c_0$ and $|\Lambda_{j_0}|<\bar{N}$, then update $\Lambda_{j_0}=\Lambda_{j_0}\cup x$ and $\Lambda_j^2=\Lambda_j^2\backslash x$.

$\bullet$ Step 5. Repeat Step 4 until any sample $x$ in the set $\bigcup\limits_{j=1}^l\Lambda_j^2$ satisfies $q_{\Lambda_{k}\cup x}\leq c_0$ for $k=l+1,\cdots,m$.

$\bullet$ Step 6: Update $\Lambda_j=\Lambda_j^1\cup\Lambda_j^2$ for $j=1,\cdots,l$.

The obtained subsets $\Lambda_1,\dots,\Lambda_m$ are the results after sample number equalization and are applied to DKI. 

\begin{figure*}[t]
    \centering
    \subfigcapskip=-2pt
    \subfigure[Samples]{\includegraphics[width=3.1cm, height=3.1cm]{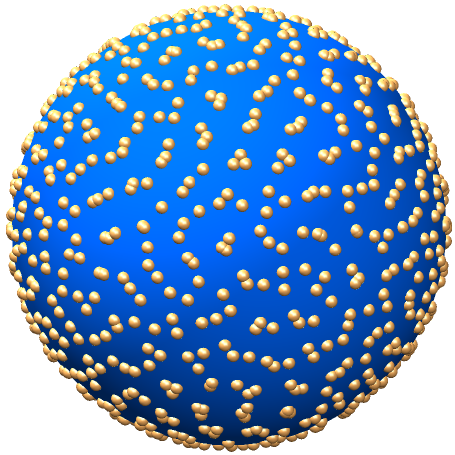}}
    \subfigure[Random]{\includegraphics[width=3.1cm, height=3.1cm]{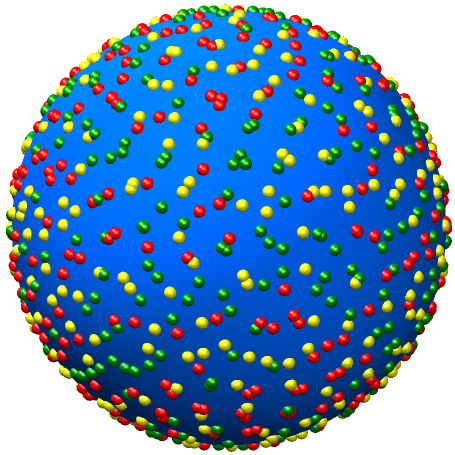}}
    \subfigure[$\tau$-uniform]{\includegraphics[width=3.1cm, height=3.1cm]{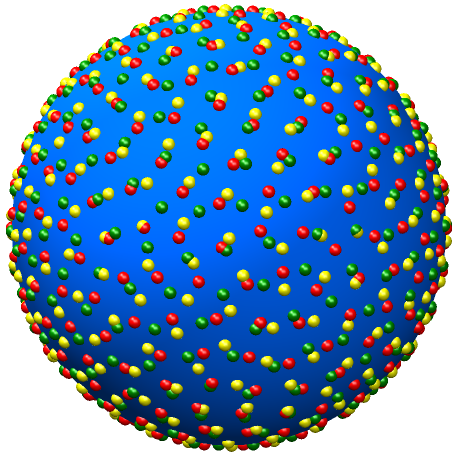}}
    \subfigure[SAJ]{{\includegraphics[width=3.1cm, height=3.1cm]{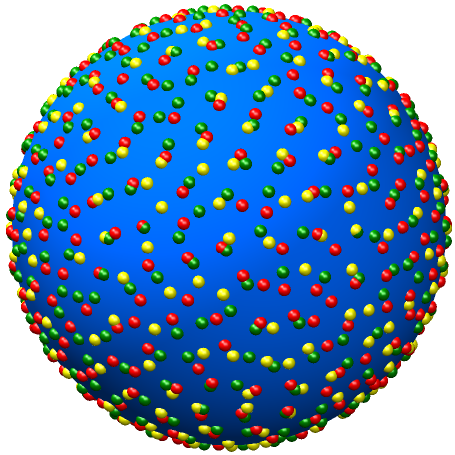}}}
	\caption{Illustration for the divisions of samples. Womersley's symmetric spherical 25 designs are rotated using 3 rotation matrices to generate the samples}\label{DataDivision3}
\end{figure*}

Figures \ref{DataDivision3}   presents intuitive comparisons among different divisions. Here, ``Random'' denotes to randomly select the points from the ``samples''. $\tau$-uniform denotes to divide the ``samples'' into three $25$-designs. ``SAJ'' denotes the proposed select-and-judge  strategy.
It can be found in these figures that SAJ is easy to yield quasi-uniform data on spheres.

\begin{figure*}[t]
    \centering
    \subfigure[Number of samples and RMSE]{\includegraphics[width=6cm,height=4.5cm]{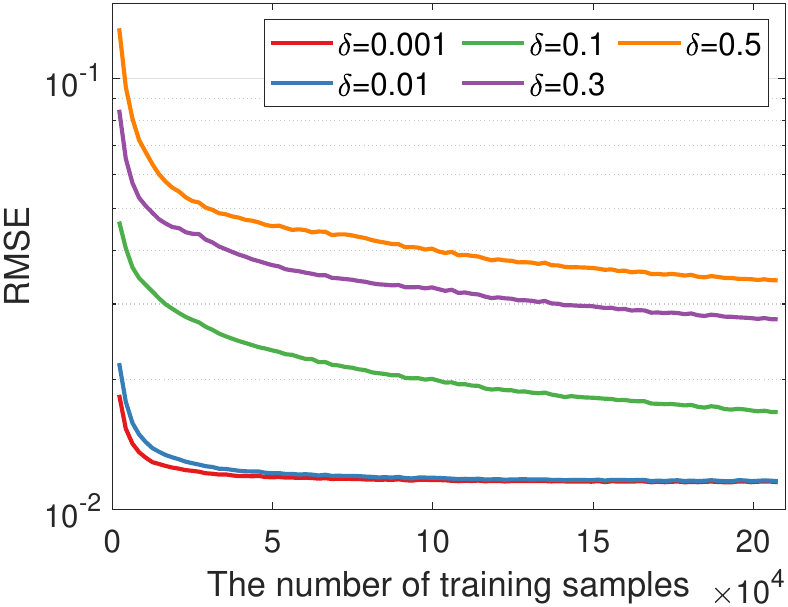}}
    \subfigure[Number of local machines and RMSE]{\includegraphics[width=6cm,height=4.5cm]{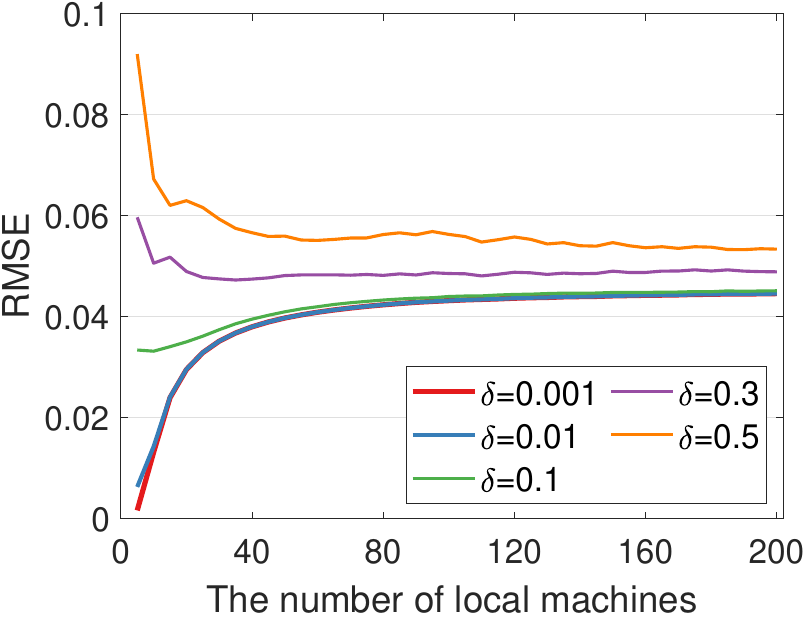}}
    \caption{Numerical results for 50-dimensional data}\label{sim2_RMSE_increase_sam}
\end{figure*}

\section{Appendix B: Simulations for high-dimensional data}
The setting of the simulation is similar as that for  $3$-dimensional data. The training samples' inputs $\{{x}_i\}_{i=1}^N$ with $N=1000$ and the testing samples' inputs $\{x_j'\}_{j=1}^{N'}$ and $N'=10000$ are i.i.d. drawn according to the uniform distribution on the $\mathbb S^{50}$  and  $\phi_2(x_1, x_2)=\exp\left\{-\|x_1-x_2\|_2^2/(2\sigma^2)\right\}$ is selected as the kernel function with $\sigma$ (equally spaced logarithmic)    chosen from 20 values drawn in  
$[0.1, 100]$.  The number of local machines $m$ varies from the set $\{2,4,\cdots,200\}$. The numerical results can be found in Figure \ref{sim2_RMSE_increase_sam}, where  similar trends as those of the $3$-dimensional data are exhibited.